\documentclass[12pt,reqno]{amsart} 
\usepackage{amssymb}
\usepackage{amsmath}
\usepackage{comment}
\usepackage{amsthm}
\usepackage{caption}
\usepackage{mathtools}
\usepackage{fullpage}
\allowdisplaybreaks
\newtheorem{theorem}{Theorem}
\newtheorem{proposition}[theorem]{Proposition}
\newtheorem{lemma}[theorem]{Lemma}
\newcommand{\EndProof}{\hfill $\square$}

\usepackage[colorlinks=true, linkcolor=blue,
 citecolor=blue, urlcolor=blue]{hyperref}

\begin{document} 
\numberwithin{equation}{section}

\title{The Fourier Coefficients of a Metaplectic Eisenstein Distribution on the Double Cover of SL$(3)$ over $\mathbb{Q}$.} 
\author{Edmund Karasiewicz}
\email{ekarasie@ucsc.edu}
\address{Edmund Karasiewicz: Department of Mathematics, University of California Santa Cruz, Santa Cruz CA 95064, USA }
\begin{abstract} We compute the Fourier coefficients of a minimal parabolic Eisenstein distribution on the double cover of SL$(3)$ over $\mathbb{Q}$. Two key aspects of the paper are an explicit formula for the constant term, and formulas for the Fourier coefficients at the ramified place $p=2$. Additionally, the unramified non-degenerate Fourier coefficients of this Eisenstein distribution fit into the combinatorial description provided by Brubaker-Bump-Friedberg-Hoffstein \cite{BBFH07}. 

\end{abstract}
\maketitle


\section{Introduction} 

The study of metaplectic Eisenstein series goes back at least to Maass \cite{M37} who computed the Fourier expansion of half-integral weight Eisenstein series and found quadratic Dirichlet series among the non-degenerate coefficients. Kubota \cite{K69}, inspired by the works of Hecke \cite{H23}, Selberg \cite{S62}, and Weil \cite{W64}, took the next step by computing the Fourier coefficients of an Eisenstein series on the $n$-fold cover of GL$(2)$ over a number field containing the $2n$-th roots of unity. Using the Fourier expansion Kubota was able to study the automorphic residues of the metaplectic Eisenstein series. In the case of the double cover of GL$(2)$ the Jacobi $\theta$-function is such a residue. Subsequent works continued to include the hypothesis that the base field contains the $2n$-th roots of unity. (The construction of $n$-fold covering groups requires the base field to contain the $n$-th roots of unity.) Kazhdan-Patterson \cite{KP84} developed a general theory of automorphic forms on $n$-fold covers of GL$(r+1)$. Their work includes a computation of the Fourier coefficients of metaplectic Eisenstein series and a study of their residues. Brubaker-Bump-Friedberg-Hoffstein \cite{BBFH07} and Brubaker-Bump-Friedberg \cite{BBF11} provided a combinatorial description of the non-degenerate Fourier coefficients of Eisenstein series on an $n$-fold cover of GL$(r+1)$. However, the double cover exists over $\mathbb{Q}$ and thus it is natural to try to extend the results of the above authors to this case. 

This paper represents a modest step toward bridging the gap between our understanding of Eisenstein series on $n$-fold covers over number fields containing the $2n$-th roots of unity and those over number fields containing only the $n$-th roots of unity. Specifically, we use the results established in \cite{K18} to derive explicit formulas for the Fourier coefficients of a minimal parabolic Eisenstein series on the double cover of SL$(3)$ over $\mathbb{Q}$, as opposed to a number field containing $\mathbb{Q}(i)$. These explicit formulas provide a means to study the poles of the Eisenstein series and its residual representations, a point we will return to momentarily. 

We will now highlight some of the notable features of this work. Theorem \ref{theorem:Constant} includes an explicit formula for the constant term, which can be used to determine the poles of the Eisenstein series. Proposition \ref{metawhitwell} contains the formula for the non-degenerate Fourier coefficients and the results of Subsection \ref{sec:BBFHform} express the unramified $p$-parts in the style of Brubaker-Bump-Friedberg-Hoffstein \cite{BBFH07}. Subsection \ref{ssec:Sigma4} contains formulas for the ramified parts $(p=2)$ of the non-degenerate Fourier coefficients. It would be interesting to know if these ramified coefficients also admit a combinatorial description in the style of \cite{BBFH07}. 

We would also like to mention that this computation is performed in the context of automorphic distributions. This perspective should be advantageous for studying certain Archimedean integrals, as demonstrated by Miller-Schmid \cite{MS12}. Specifically, Rankin-Selberg integrals for the exterior square L-function on GL$(n)$ were provided by Jacquet-Shalika \cite{JS90} and Bump-Friedberg \cite{BF90}, but the Archimedean integrals were not evaluated. However, Miller-Schmid \cite{MS12} computed the Archimedean integrals of the exterior square L-function on GL$(n)$ by revisiting the construction of Jacquet-Shalika in the context of automorphic distributions. A similar situation persists in the case of the symmetric square L-function. A Rankin-Selberg construction is known, but the computation of the Archimedean integrals remains open. Thus a natural test for the theory of automorphic distributions, as it relates to Archimedean integrals, is the computation of the Archimedean integrals of the symmetric square L-function on GL$(n)$. When $n=3$, the formulas for the Fourier coefficients of the minimal parabolic Eisenstein distribution contained in this paper provide a means to study the distributional analog of the $\theta$-function (i.e. a residual representation) used in the Rankin-Selberg construction of the symmetric square L-function due to Patterson and Piatetski-Shapiro \cite{PPS89}.


We will now provide a brief description of the contents of this paper. Section \ref{notation} contains notation and basic computations. Subsection \ref{ssec:2-cocycle} introduces the Banks-Levy-Sepanski 2-cocycle \cite{BLS99} and collects some computations involving this 2-cocycle. Subsection \ref{splittingreview} reviews some results from \cite{K18} about an arithmetic splitting function $s$. These results provide the foundation for all of the subsequent computations.


Section \ref{sec:ExpSums} isolates the computations of the exponential sums associated with the big Bruhat cell appearing in the Fourier coefficients of the metaplectic Eisenstein distribution. The symmetries of $s$ studied in Subsection \ref{sec:ExpSumPrelim} induce enough symmetries of the exponential sums to reduce the general computation to more manageable special cases, which are treated in subsections \ref{ssec:Sigma2} and \ref{ssec:Sigma4}. The principal reduction is described in Proposition \ref{exptwistmult}. Another significant consequence of Proposition \ref{exptwistmult} is that the Dirichlet series appearing in the non-degenerate Fourier coefficients, which do not possess an Euler product, can be reconstructed from their $p$-parts. Similar twisted multiplicativity appears in the work of Brubaker-Bump-Friedberg-Hoffstein \cite{BBFH07}. 

Section \ref{sec:MetaEisenstein} contains the computation of the Fourier coefficients of the minimal parabolic Eisenstein distribution on the double cover of SL$(3)$ over $\mathbb{Q}$. An outline of the computation follows. We begin, as is typical with computations of Fourier coefficients of Eisenstein series, by breaking up the Eisenstein distribution into six pieces using the Bruhat decomposition. The computation of the Fourier coefficients of each of the six pieces occupies one of the six subsections \ref{ssec:tildeB}-\ref{ssec:tildeNwlB}. Specifically, the formulas for the Fourier coefficients can be found in Proposition \ref{metawhitid}, Proposition \ref{metawhitalpha1}, Proposition \ref{metawhitalpha2}, Proposition \ref{metawhitalpha12}, Proposition \ref{metawhitalpha21}, and Proposition \ref{metawhitwell}. For each of the six pieces, we compute the Fourier coefficients using the method of unfolding which reduces the computation to the determination of certain exponential sums. The most complicated exponential sums appear in the case of the big Bruhat cell, but these are exactly the exponential sums studied in Section \ref{sec:ExpSums}; the exponential sums appearing in the degenerate coefficients are related to Fourier coefficients of Eisenstein series on the double cover of SL$(2)$ over $\mathbb{Q}$ and are handled directly. Finally, the formula for the constant term of the Eisenstein distribution appears in Theorem \ref{theorem:Constant} in Subsection \ref{ssec:Constant}. This formula is derived by specializing the computations of the Fourier coefficients in subsections \ref{ssec:tildeB}-\ref{ssec:tildeNwlB} to the case of the constant term.


\section{Notation} \label{notation}\label{sec:Notation}\label{sec:Comp}
\subsection{SL($3,\mathbb{R}$) and $\widetilde{\text{SL}}(3,\mathbb{R})$}\label{ssec:LieGroups}

This section contains the notation and basic computations that will be used throughout this paper. 


Let $\widetilde{\text{SL}}(3,\mathbb{R})$ be the nontrivial topological double cover of $\text{SL}(3,\mathbb{R})$. As a set $\widetilde{\text{SL}}(3,\mathbb{R})\cong \text{SL}(3,\mathbb{R})\times\{\pm1\}$. 
The Banks-Levy-Sepanski 2-cocycle $\sigma:\text{SL}(3,\mathbb{R})\times\text{SL}(3,\mathbb{R})\rightarrow \{\pm1\}$ \cite{BLS99}, recalled in Subsection \ref{ssec:2-cocycle}, defines the group multiplication on $\widetilde{\text{SL}}(3,\mathbb{R})$ as follows:
\begin{equation}\label{metamult}(g_1,\epsilon_1)(g_2,\epsilon_2) = (g_1g_2,\epsilon_1\epsilon_2\sigma(g_1,g_2)).
\end{equation}

The following list establishes notation for some subgroups of SL$(3,\mathbb{R})$ and $\widetilde{\text{SL}}(3,\mathbb{R})$:
\[\def\arraystretch{1.2}
\begin{array}{lclclcl}
G&=&\text{SL}(3,\mathbb{R}) &, & \widetilde{G}&=&\widetilde{\text{SL}}(3,\mathbb{R}),\\
B&=&\left\{\left(\begin{smallmatrix}
a & b & c\\
0 & e & f\\
0 & 0 & \frac{1}{ae}\end{smallmatrix}\right)| a,e\neq0\right\} &, &
\widetilde{B}&=&\left\{\left(\left(\begin{smallmatrix}
a & b & c\\
0 & e & f\\
0 & 0 & \frac{1}{ae}\end{smallmatrix}\right),\pm1\right)| a,e\neq0\right\},\\
N&=&\left\{\left(\begin{smallmatrix}
1 & x & z\\
0 & 1 & y\\
0 & 0 & 1\end{smallmatrix}\right)\right\} &, &
\widetilde{N}&=&\left\{\left(\left(\begin{smallmatrix}
1 & x & z\\
0 & 1 & y\\
0 & 0 & 1\end{smallmatrix}\right),1\right)\right\},\\
T&=&\left\{\left(\begin{smallmatrix}
t_1 & 0 & 0\\
0 & t_2 & 0\\
0 & 0 & \frac{1}{t_1t_2}\end{smallmatrix}\right)|t_i\in\mathbb{R}^{\times}\right\}&, &
\widetilde{T}&=&\left\{\left(\left(\begin{smallmatrix}
t_1 & 0 & 0\\
0 & t_2 & 0\\
0 & 0 & \frac{1}{t_1t_2}\end{smallmatrix}\right),\pm1\right)|t_i\in\mathbb{R}^{\times}\right\},\\
A&=&\left\{\left(\begin{smallmatrix}
a & 0 & 0\\
0 & b & 0\\
0 & 0 & \frac{1}{ab}\end{smallmatrix}\right)|a,b>0\right\}&, &
\widetilde{A}&=&\left\{\left(\left(\begin{smallmatrix}
a & 0 & 0\\
0 & b & 0\\
0 & 0 & \frac{1}{ab}\end{smallmatrix}\right),1\right)|a,b>0\right\},\\
M&=&\left\{\left( \begin{smallmatrix}
\epsilon_1 & 0 & 0\\
0 & \epsilon_2 & 0\\
0 & 0 & \epsilon_1\epsilon_2 \end{smallmatrix} \right)|\epsilon_1,\epsilon_2=\pm1\right\}&, &
\widetilde{M}&=&\left\{\left( \left(\begin{smallmatrix}
\epsilon_1 & 0 & 0\\
0 & \epsilon_2 & 0\\
0 & 0 & \epsilon_1\epsilon_2 \end{smallmatrix} \right),\pm1\right)|\epsilon_1,\epsilon_2=\pm1\right\},\\
K&=&\text{SO}(3)&,&\widetilde{K}&=&\text{Spin}(3),
\end{array}\]
\begin{equation}\label{gamma14+infty}\begin{array}{lll}
\Gamma=\Gamma_{1}(4)= 
\{\gamma\in \text{SL}(3,\mathbb{Z})|\gamma \equiv \left(\begin{smallmatrix}
1&*&*\\
0&1&*\\
0&0&1\end{smallmatrix}
\right)\text{(mod }4)\} &, & \Gamma_{\infty}=N\cap\text{SL}(3,\mathbb{Z}).
\end{array}\end{equation}

To simplify notation $\left(\left(\begin{smallmatrix}
a & 0 & 0\\
0 & b & 0\\
0 & 0 & \frac{1}{ab}\end{smallmatrix}\right),1\right)$ may be written 
$\left(\begin{smallmatrix}
a & 0 & 0\\
0 & b & 0\\
0 & 0 & \frac{1}{ab}\end{smallmatrix}\right),$ and 
$\left(\left(\begin{smallmatrix}
1 & x & z\\
0 & 1 & y\\
0 & 0 & 1\end{smallmatrix}\right),1\right)$ may be written
$\left(\begin{smallmatrix}
1 & x & z\\
0 & 1 & y\\
0 & 0 & 1\end{smallmatrix}\right)$. Furthermore, we will write $
\text{t}(a,b,c) = \left(\begin{smallmatrix}
a & 0 & 0\\
0 & b & 0 \\
0 & 0 & c
\end{smallmatrix}\right)$, and 
$n(x,y,z) = \left(\begin{smallmatrix}
1 & x & z\\
0 & 1 & y\\
0 & 0 & 1\end{smallmatrix}
\right)
$.

Now we list particular representatives of elements of the Weyl group of SL$(3,\mathbb{R})$:
\[\begin{array}{lllllll}
&w_{\alpha_1}&=&\left(\begin{smallmatrix}
0 & -1 & 0\\
1 & 0 & 0\\
0 & 0 & 1\end{smallmatrix}\right),&
w_{\alpha_2}&=&\left(\begin{smallmatrix}
1 & 0 & 0\\
0 & 0 & -1\\
0 & 1 & 0\end{smallmatrix}\right),\\
&w_{\alpha_1}w_{\alpha_2}&=&\left(\begin{smallmatrix}
0 & 0 & 1\\
1 & 0 & 0\\
0 & 1 & 0\end{smallmatrix}\right),&
w_{\alpha_2}w_{\alpha_1}&=&\left(\begin{smallmatrix}
0 & -1 & 0\\
0 & 0 & -1\\
1 & 0 & 0\end{smallmatrix}\right),\\
\text{and}&
w_{\ell}&=&w_{\alpha_1}w_{\alpha_2}w_{\alpha_1}&=w_{\alpha_2}w_{\alpha_1}w_{\alpha_2}&=&\left(\begin{smallmatrix}
0 & 0 & 1\\
0 & -1 & 0\\
1 & 0 & 0\end{smallmatrix}\right).
\end{array}\]
These representatives of the Weyl group are used in the formula for the 2-cocycle defined in \cite{BLS99}. We recall this formula in Subsection \ref{ssec:2-cocycle}.

If $g\in\text{SL}(3,\mathbb{R})$, let $^\intercal g$ denote the transpose of $g$. If $Q$ is a subgroup of $B$, then $Q_{-}$ and $Q^{op}$ will denote the set $\{^\intercal q| q\in Q\}$.

Let $\frak{sl}(3,\mathbb{R})$ be the real Lie algebra of SL$(3,\mathbb{R})$ and let $\frak{a}$ be the subalgebra of diagonal matrices. For $X\in \frak{sl}(3,\mathbb{R})$, let exp$(X) = \sum_{n=0}^{\infty}\frac{X^n}{n!}$ be the exponential map and let log denote its inverse on $\frak{a}$. The map $K\times A \times N\rightarrow \text{SL}(3,\mathbb{R})$ defined by multiplication is a diffeomorphism. Define the maps $\kappa:\text{SL}(3,\mathbb{R})\rightarrow K,\,H:\text{SL}(3,\mathbb{R})\rightarrow \frak{a}$, and $\nu:\text{SL}(3,\mathbb{R})\rightarrow N$ such that $g\mapsto (\kappa(g),\text{exp}(H(g)),\nu(g))$ is the inverse of the map $(k,a,n)\mapsto kan$.

\subsection{The $A_2$ Root System} \label{ssec:A2}
This subsection establishes notation related to the root system of SL$(3,\mathbb{R})$. Let $X(T)=\mathbb{Z}e_1+\mathbb{Z}e_2+\mathbb{Z}e_3/\mathbb{Z}(e_1+e_2+e_3)$ be the group of rational characters of $T$ (written additively), where $e_i(\text{t}(t_1,t_2,t_3))=t_i$. The Weyl group $W$ acts on $X(T)$ through its action on $T$ by conjugation; this action extends to an action on the vector space $X(T)\otimes \mathbb{C}$. 

Let $\Phi=\Phi(\text{SL}(3,\mathbb{R}),T)=\Phi^{+}\cup\Phi^{-}$, where $\Phi^{+}=\{(e_1-e_2),(e_2-e_3),(e_1-e_3)\}$ is the  set of positive roots of SL$(3,\mathbb{R})$ with respect to $B$ and $\Phi^{-}=-\Phi^{+}$ is the set of negative roots. To each root $\alpha=\pm(e_i-e_j)\in \Phi$ there is a canonically defined element $h_{\alpha}=\pm(e_{ii}-e_{jj})\in\frak{a}$, where $e_{ij}$ is the $3\times3$ matrix with $1$ in the $(i,j)$-position and $0$ elsewhere; details can be found in \cite{H72}. The elements of $X(T)\otimes \mathbb{C}$ pair with these matrices by the formula $\langle a_1e_1+a_2e_2+a_3e_3,\pm(e_{ii}-e_{jj})\rangle=\pm(a_{i}-a_j)$. This pairing shows that $X(T)\otimes \mathbb{C}$ is isomorphic to $\frak{a}^{\prime}_{\mathbb{C}}$, the space of complex valued linear functionals of $\frak{a}$. It will be convenient to embed $X(T)\otimes \mathbb{C}$ into $\mathbb{C}^3$ via the map $a_1e_1+a_2e_2+a_3e_3\mapsto (\frac{2}{3}a_1-\frac{1}{3}a_2-\frac{1}{3}a_3, -\frac{1}{3}a_1+\frac{2}{3}a_2-\frac{1}{3}a_3,-\frac{1}{3}a_1-\frac{1}{3}a_2+\frac{2}{3}a_3)$. In this case the pairing with $\frak{a}$ is given by $\langle (\lambda_1,\lambda_2,\lambda_3),\pm(e_{ii}-e_{jj})\rangle=\pm(\lambda_{i}-\lambda_j)$.

\subsection{Pl\"{u}cker Coordinates} \label{sec:SL3R} \label{sec:PluckerSym}

Given $\left( \begin{smallmatrix}
a & b & c \\
d & e & f \\
g & h & i \end{smallmatrix} \right)\in $ SL$(3,\mathbb{R})$, define six parameters, called Pl\"{u}cker coordinates, as follows:
\begin{equation*}
A_1^\prime = -g, \hspace{.5cm} B_1^\prime = -h, \hspace{.5cm} C_1^\prime = -i,\hspace{.5cm}
A_2^\prime = -(dh-eg), \hspace{.5cm} B_2^\prime = (di-fg), \hspace{.5cm} C_2^\prime = -(ei-fh).
\end{equation*}
\begin{theorem} \label{theorem:PluckerCoords}(\cite[Ch 5]{BGL3}) The map taking $\left( \begin{smallmatrix}
a & b & c \\
d & e & f \\
g & h & i \end{smallmatrix} \right)\in $ SL$(3,\mathbb{R})$ to $(A_1^\prime,B_1^\prime,C_1^\prime,A_2^\prime,B_2^\prime,C_2^\prime)$ defines a bijection between the coset space $N\backslash$SL$(3,\mathbb{\mathbb{R}})$ and the set of all \newline$(A_1^\prime,B_1^\prime,C_1^\prime,A_2^\prime,B_2^\prime,C_2^\prime)\in \mathbb{R}^6$ such that: $A_1^\prime C_2^\prime+B_1^\prime B_2^\prime+C_1^\prime A_2^\prime= 0$, not all of $A_1^\prime,B_1^\prime,C_1^\prime$ equal $0$, and not all of $A_2^\prime,B_2^\prime,C_2^\prime$ equal $0$.  Furthermore, a coset in $N\backslash$SL$(3,\mathbb{\mathbb{R}})$ contains an element of SL$(3,\mathbb{Z})$ if and only if $A_1^\prime,B_1^\prime,C_1^\prime$ are coprime integers and $A_2^\prime,B_2^\prime,C_2^\prime$ are coprime integers.\end{theorem}

Versions of this result hold for other congruence subgroups. Let
$A_1^\prime = 4A_1$,  $A_2^\prime = 4A_2$,
$B_1^\prime = 4B_1$,  $B_2^\prime =4B_2$,
$C_1^\prime = C_1$,  $C_2^\prime = C_2$.
The coset space 
 $\Gamma_{\infty}\backslash \Gamma_1(4)$ can be identified with
\begin{equation}\label{cosetparam2}\left\{\begin{array}{c}
(4A_1,4B_1,C_1,4A_2,4B_2,C_2)\in \mathbb{Z}^6 | A_1C_2+4B_1B_2+C_1A_2 = 0,\\
(A_i,B_i,C_i) = 1,\,C_j\equiv -1\text{ (mod } 4)
\end{array}\right\}
\begin{array}{l}
\text{}\\
.\end{array}
\end{equation}

\noindent\textbf{Coset Representatives:}
The following table lists coset representatives of $\Gamma_{\infty}\backslash \Gamma$ following Bump \cite{BGL3}.

\begin{center}
\begin{tabular}{|c|c|c|}\hline
 Cell & Constraints & $\Gamma_{\infty} \backslash \Gamma$ Representative\\ \hline
 $B$ & $C_1,C_2\neq 0$ & $\left( \begin{smallmatrix}
\frac{-1}{C_2} &  & \\
 & \frac{C_2}{C_1} & \\
 &  & -C_1\end{smallmatrix} \right)$ \\ \hline
 $Nw_{\alpha_1}B$ & $\begin{array}{c}A_1,B_1,A_2 = 0,\\ C_1,B_2\neq 0\end{array}$ & $\left( \begin{smallmatrix}
 & -1 & \\
-1 &  & \\
 &  &  -1\end{smallmatrix} \right)
\left( \begin{smallmatrix}
\frac{4B_2}{C_1} & \frac{-C_2}{C_1} & \\
 & \frac{1}{4B_2} & \\
 &  & C_1 \end{smallmatrix} \right)$\\ \hline
 $Nw_{\alpha_2}B$ & $\begin{array}{c}A_1,A_2,B_2 = 0,\\ B_1,C_2\neq 0\end{array}$& $\left( \begin{smallmatrix}
 -1 &  & \\
 &  & -1\\
 & -1 &  \end{smallmatrix} \right)
\left( \begin{smallmatrix}
\frac{1}{C_2} &  & \\
 & 4B_1 & C_1\\
 &  & \frac{C_2}{4B_1} \end{smallmatrix} \right)$\\ \hline
 $Nw_{\alpha_1}w_{\alpha_2}B$ & $\begin{array}{c}A_1= 0,\\ B_1,A_2\neq 0\end{array}$ & $\left( \begin{smallmatrix}
  & & 1\\
1 &  & \\
& 1 &  \end{smallmatrix} \right)
\left( \begin{smallmatrix}
\frac{A_2}{B_1} &  & \frac{C_2}{4B_1}\\
 & -4B_1 & \frac{4B_1 B_2}{A_2}\\
 &  & \frac{-1}{4A_2} \end{smallmatrix} \right)$\\ \hline
 $Nw_{\alpha_2}w_{\alpha_1}B$ & $\begin{array}{c}A_2= 0,\\ A_1,B_2\neq 0\end{array}$& $\left( \begin{smallmatrix}
  & 1 & \\
 &  & 1\\
1 &  &  \end{smallmatrix} \right)
\left( \begin{smallmatrix}
-4A_1 & -4B_1 & -C_1\\
 & \frac{-1}{4B_2} & \\
 &  & \frac{B_2}{A_1} \end{smallmatrix} \right)$\\ \hline
$Nw_{\ell}B$ & $A_1,A_2\neq 0$ &  $\left( \begin{smallmatrix}
 &  & -1\\
 & -1 & \\
-1 & & \end{smallmatrix} \right)
\left( \begin{smallmatrix}
4A_1 & 4B_1 & C_1\\
 & \frac{A_2}{A_1} & \frac{-B_2}{A_1}\\
 &  & \frac{1}{4A_2} \end{smallmatrix} \right)$\\ \hline
 \end{tabular}
 \end{center}
 \begingroup
 \captionof{table}{$\Gamma_{\infty}\backslash\Gamma$ representatives}\label{CosetReps}
 \endgroup
 

The next proposition collects some symmetries satisfied by the Pl\"{u}cker coordinates.
\begin{proposition} \label{prop:PluckerSym}
Let $g\in \text{SL}(3,\mathbb{R})$ with Pl\"{u}cker coordinates $(4A_1,4B_1,C_1,4A_2,4B_2,C_2)$, \newline$n=n(x,y,z)\in N$, $S_2=\mathrm{t}(1,-1,1)$, and $S_3 = \mathrm{t}(1,1,-1)$. Then:
\begin{enumerate}
\item The matrix $ngn^{-1}$ has Pl\"{u}cker coordinates \[(4A_1,4B_1-4A_1x,C_1-4B_1y+4A_1(xy-z),4A_2,4B_2+4A_2y,C_2+4B_2x+4A_2z).\]
\item The matrix $S_3gS_3^{-1}$ has Pl\"{u}cker coordinates $(-4A_1,-4B_1,C_1,-4A_2,4B_2,C_2)$.
\item The matrix $S_2gS_2^{-1}$ has Pl\"{u}cker coordinates $(4A_1,-4B_1,C_1,4A_2,-4B_2,C_2)$.
\item The matrix $w_{\ell}\,^\intercal g^{-1}w_{\ell}$ has Pl\"{u}cker coordinates $(4A_2,-4B_2,C_2,4A_1,-4B_1,C_1)$.
\item Let $g\in\Gamma_{1}(4)$. If $D$ divides $(A_1,A_2)$, $D_1=(D,B_1)$, $D=D_1D_2$, and \newline $T=\text{t}(1,D_2^{-1},D^{-1})$ then $TgT^{-1}\in\text{SL}(3,\mathbb{Z})$ has Pl\"{u}cker coordinates 
\[(4A_1/D,4B_1/D_1,C_1,4A_2/D,(4B_2)/D_2,C_2).\]
Furthermore, $TgT^{-1}\in\Gamma_1(4)$ if and only if $D_2$ divides $B_2$. 
\end{enumerate}
\end{proposition}
\noindent The proof is straightforward matrix algebra and will be omitted.

Let 
\begin{equation} \label{BoldSA1A2}
\mathbb{S}(A_1,A_2)=\{\gamma\in\Gamma_1(4)| \gamma \text{ has Pl\"{u}cker coordinates of the form }(4A_1,*,*,4A_2,*,*)\},\end{equation} 
and let
\begin{multline}\label{SA1A2}S(A_1,A_2) = \{(4A_1,4B_1,C_1,4A_2,4B_2,C_2)\in \mathbb{Z}^6 |
A_1C_2+4B_1B_2+C_1A_2 = 0,\\
(A_i,B_i,C_i) = 1,\text{ } C_j\equiv -1\text{(mod } 4), \text{ } \frac{B_1}{A_1},\frac{B_2}{A_2},\frac{C_2}{4A_2}\in[0,1)\}.
\end{multline}
The maps of Proposition \ref{prop:PluckerSym} induce maps on the double coset space $\Gamma_\infty\backslash \mathbb{S}(A_1,\mu A_2)/\Gamma_\infty\cong S(A_1,A_2)$. ((1) in Proposition \ref{prop:PluckerSym} implies the bijection.) The next proposition describes another important property of these double coset spaces. We will refer to this property as multiplicativity of the double coset spaces.

\begin{proposition}[\cite{K18}] \label{mult} Let $A_1,\alpha_1> 0$, $A_2,\alpha_2\neq0,$ suppose that $(A_1A_2,\alpha_1\alpha_2) = 1$, $A_1,A_2$ are odd, and suppose that $A_1\alpha_1+A_2\alpha_2\equiv 0$ (mod 4). Let $\mu=(\frac{-1}{-A_1A_2})$. Then 
\[\Gamma_\infty\backslash \mathbb{S}(A_1\alpha_1,A_2\alpha_2)/\Gamma_\infty \cong \Gamma_\infty\backslash \mathbb{S}(A_1,\mu A_2)/\Gamma_\infty \times\Gamma_\infty\backslash \mathbb{S}(\alpha_1,-\mu\alpha_2)/\Gamma_\infty.\]

The bijection is induced by the map
\begin{align}\nonumber(4A_1\alpha_1,4B_1&,C_1,4A_2\alpha_2,4B_2,C_2)\stackrel{\varphi}{\mapsto}\\
\nonumber((4A_1,4B_1,\mathcal{C}_1,\mu4 A_2,4B_2,\gamma &C_2),\\
(4\alpha_1,&4B_1,\Big(\frac{-1}{A_2}\Big) A_2C_1,-\mu4\alpha_2,-\Big(\frac{-1}{A_2}\Big)\mu4B_2,-\mu \Big(\frac{-1}{A_2}\Big)A_1C_2)),
\end{align}
where:
\begin{enumerate}
\item$\mathcal{C}_1=\frac{-A_1\gamma C_2-4B_1B_2}{\mu A_2}$;
\item$\gamma$ is the smallest positive integer such that $\gamma\equiv 1$ (mod 4) and $\gamma\equiv \alpha_1$ (mod $A_2$).
\end{enumerate}
\end{proposition}

\subsection{Exponential Sums}\label{ssec:ExpSums}
This section collects some basic identities involving Gauss and Ramanujan sums. The proofs are elementary and will be omitted. The definition and basic properties of the Kronecker symbol $(\frac{\cdot}{\cdot})$ can be found in \cite{K18} (or \cite{IK04}).

When $d$ divides $n$, let
\begin{equation}\label{GSum}g(d,m,n)= \sum_{x\in (\mathbb{Z}/n\mathbb{Z})^\times}(\frac{x}{d})e^{2\pi i(\frac{mx}{n})}.
\end{equation}
When $n=2^k$, $k\geq 3$, and $\epsilon\equiv \pm1\,(\text{mod }4)$, let 
\begin{equation}\label{GSumEven}g_{\epsilon}(2^{i},m,2^{k}) = \sum_{\begin{smallmatrix}
x\in (\mathbb{Z}/2^k\mathbb{Z})^\times\\
x\equiv \epsilon\, (\text{mod }4)
\end{smallmatrix}}\big(\frac{x}{2^i}\big)e^{2\pi i(m\frac{x}{2^k})}.
\end{equation}
Note that if $i$ is even then $g_{\epsilon}(2^{i},m,2^{k})$ is also well defined for $k=2$.

In what follows, if $P$ is a statement then we define $\delta_{P}$ to be $1$ if $P$ is true; $0$ if $P$ is false. If $a,m\in \mathbb{Z}$, we write $\bar{a}$ or $(a)^{-1}_{\text{mod }m}$ for a multiplicative inverse of $a$ modulo $m$.

\begin{lemma} \label{lemma:gauss} Let $p$ be an odd prime, let $a,b,m\in\mathbb{Z}$, and let $j,k,\ell \in\mathbb{Z}_{\geq 0}$.
\begin{enumerate} 
\item \label{Gauss1}If $c\in(\mathbb{Z}/d\mathbb{Z})^{\times}$, then $g(d,cm,n) = \left(\frac{c}{d}\right)g(d,m,n)$.
\item \label{Gauss2}If $k\neq0$, then
\begin{align*} g(p^1,p^j,p^k) &=
\begin{cases} 
      0, & k\neq j+1; \\
       p^{k-1}g(p^1,p^0,p^1), & k=j+1, \\ 
   \end{cases}\\
g(p^0,p^j,p^k)&=
\begin{cases} 
      \phi(p^k), & k-j\leq 0; \\
       -p^{k-1},& k-j=1; \\
      0, & k-j\geq 2. 
   \end{cases}
\end{align*}
\item\label{Gauss3} $p^\ell g(p^i,\pm p^j,p^k) = g(p^i,\pm p^{j+\ell},p^{k+\ell})$.
\item\label{Gauss3even} $2^\ell g_{\epsilon}(2^i,\pm 2^j,2^k) = g_{\epsilon}(2^i,\pm 2^{j+\ell},2^{k+\ell})$.
\item\label{Gauss4} Suppose that $p$ does not divide $a$ and $j>0$, then 
\begin{equation}\label{AffineGaussSum}
\sum_{x\in (\mathbb{Z}/p^j\mathbb{Z})^\times}(\frac{ax+b}{p})^k e^{2\pi i\frac{mx}{p^j}}=e^{2\pi i\frac{-m\bar{a}b}{p^j}}(\frac{a}{p})^kg(p^k,m,p^j)-(\frac{b}{p})^k\sum_{\ell=0}^{p^{j-1}-1} e^{2\pi i\frac{m\ell}{p^{j-1}}}.
\end{equation}
\item \begin{equation}\label{Gauss5}\sum_{0\leq \ell \leq k}g(p^2,m,p^\ell)=\delta_{p^{k}|m}p^{k}.
\end{equation}
\end{enumerate}\end{lemma}

\subsection{Zeta-Functions}\label{ssec:Zeta}

This section collects some definitions and identities involving zeta-functions. Let $\phi$ be the Euler $\phi$-function and let 
\begin{equation*}K_\kappa(n;4c)=\sum_{d\in\mathbb{Z}/(4c)\mathbb{Z}
}\varepsilon_{d}^{-\kappa}(\frac{4c}{d})e^{2\pi i\frac{nd}{4c}}, \text{ where } \varepsilon_{d}=\begin{cases}
1, & \text{if } d\equiv 1\,(\text{ mod }4);\\
i, & \text{if } d\equiv -1\,(\text{ mod }4);\\
0, & \text{otherwise}.
\end{cases}\end{equation*}

The prototypical zeta-function is the Riemann zeta function
\begin{equation}\label{zeta} \zeta(s) = \prod_{
\begin{smallmatrix}
p\\
\text{prime}
\end{smallmatrix}}(1-p^{-s})^{-1}=\sum_{n>0}n^{-s},
\end{equation}
where Re$(s)>1$.
When the $2$-part is missing we will write $\zeta_2(s) = \prod_{p\neq2}(1-p^{-s})^{-1}=\zeta(s)(1-2^{-s})$.

The computation of the semi-degenerate Fourier coefficients will employ the identity
\begin{equation}\label{squarephi}
\sum_{k>0}\frac{\phi(4k^2)}{2}(2k)^{-2s}=\frac{2^{-2s}}{1-2^{-(2s-2)}}\prod_{p\neq2}\frac{1-p^{1-2s}}{1-p^{2-2s}}=2^{-2s}\frac{\zeta(2s-2)}{\zeta_2(2s-1)},
\end{equation}
and will include the appearance of the zeta-function
\begin{equation}\label{KloostermanDirichlet}a_{\epsilon,\nu}(n)=\epsilon i4^{-\nu-1}\zeta_{2}(2\nu+1)\sum_{c\in\mathbb{Z}_{>0}}c^{-\nu-1}K_{\epsilon}(-n;4c),
\end{equation}
where Re$(\nu)>1$.
Bate \cite[pg 26, 35]{B} shows that (\ref{KloostermanDirichlet}) is a Fourier coefficient of a metaplectic Eisenstein series on the double cover of $\text{SL}(2,\mathbb{R})$ and evaluates it in terms of quadratic L-functions when $n\neq0$ and $\zeta$ when $n=0$. 

The following identity is used in the calculation of the constant term. Let $n\in\mathbb{Z}_{>0}$, then

\begin{equation}\label{KloostermanConstant}
K_{\kappa}(0,4n)=\begin{cases}
0  &,\text{ if } n \text{ is not a square;}\\
(1+i^{-\kappa})\frac{\phi(4n)}{2} &,\text{ if } n \text{ is a square.}
\end{cases}
\end{equation}

\subsection{2-cocycle}\label{ssec:2-cocycle}

This section collects basic facts about the Banks-Levy-Sepanski 2-cocycle \cite{BLS99}. For $x,y\in\mathbb{R}^{\times}$ we will write $(x,y)_{\mathbb{R}}=(x,y)$ for the Hilbert symbol which is equal to $-1$ if both $x$ and $y$ are negative and $1$ otherwise. As a warning, for $x,y\in\mathbb{Z}$ we also write $(x,y)$ for the GCD of $x$ and $y$. The correct interpretation of $(x,y)$ will be clear from context.

Let $g\in $ SL$(3,\mathbb{R})$ with Pl\"{u}cker coordinates $(A_1^\prime,B_1^\prime,C_1^\prime,A_2^\prime,B_2^\prime,C_2^\prime)$. Let $X_1(g)=\text{det}(g)$, let $X_2(g)$ be the first non-zero element of the list $-A_2,B_2,-C_2$, let $X_3(g)$ be the first nonzero element of the the list $-A_1,-B_1,-C_1$, and let 
$\Delta(g)=\text{t}(X_{1}(g)/X_{2}(g),X_{2}(g)/X_{3}(g),X_{3}(g)).$

Let $\sigma$ be the 2-cocycle defined in \cite{BLS99}. The cocycle $\sigma$ can be computed as follows. If $g_1,g_2\in G$ such that $g_1=naw_{1}\ldots w_{k}n^\prime$ is the Bruhat decomposition of $g_1$, then in Section 4 of \cite{BLS99} Banks-Levy-Sepanski show that the 2-cocycle $\sigma$ satisfies the formula
\begin{equation}\label{BLSdef}\sigma(g_1,g_2) = \sigma(a,w_{1}\ldots w_{k}n^\prime g_2)\sigma(w_{1},w_{2}\ldots w_{k}n^\prime g_2)\ldots\sigma(w_{k-1},w_{k}n^\prime g_2)\sigma(w_{k},n^\prime g_2).\end{equation}
Each factor can be computed using the following rules: let $h\in G$, $a\in T$, then
\begin{align*}
\sigma(\text{t}(a_1,a_2,a_3),\text{t}(b_1,b_2,b_3))=& (a_1,b_2)(a_1,b_3)(a_2,b_3),\\
\sigma(a,h) =& \sigma(a,\Delta(h)),\\\text{ and }\hspace{1cm}
\sigma(w_\alpha,h) =& \sigma(\Delta(w_\alpha h)\Delta(h),-\Delta(h)).
\end{align*} 

The next lemma describes some simple identities involving $\sigma$.
\begin{lemma} \label{lemma:cocycle}(Banks-Levy-Sepanski \cite{BLS99})
Let $n,n_1,n_2\in N$ and let $a\in A$, then\\
$\sigma(n_1g_1,g_2n_2)=\sigma(g_1,g_2)$,
$\sigma(g_1n,g_2) = \sigma(g_1,ng_2)$,
$\sigma(n,g) = \sigma(g,n) = 1$, and
$\sigma(g,a)=1$.
\end{lemma}

Note that $\{(^{\intercal}n(x,y,z),1)\}$ is not a subgroup of $\widetilde{G}$. For this reason we define the subgroup $\widetilde{N}_{-}=(w_{\ell},1)^{-1}\widetilde{N}(w_{\ell},1)=(w_{\ell},1)\widetilde{N}(w_{\ell},1)^{-1}$.

Finally we will collect several identities that will be useful later. The proofs follow directly from the definition of $\sigma$ and are omitted.
\begin{proposition}\label{metamultwell} Let $\gamma^\prime\in \text{SL}(3,\mathbb{R})$ with Pl\"{u}cker coordinates $(A_1,B_1,C_1,A_2,B_2,C_2)$. \newline
\begin{itemize}
\item If $\gamma^\prime=\Bigg(\begin{smallmatrix}
 0 & 0 & -1\\
0 & -1 & 0\\
-1 & 0 & 0 \end{smallmatrix} \Bigg)
\Bigg(\begin{smallmatrix}
4A_1 & 4B_1 & C_1\\
 0& \frac{A_2}{A_1} & \frac{-B_2}{A_1}\\
 0& 0 & \frac{1}{4A_2} \end{smallmatrix} \Bigg)$ and $\epsilon=(A_1,-A_2)$, then
\begin{multline}(\gamma^\prime, 1)^{-1}(w_{\ell},1)=\\
\Bigg( \Bigg( \begin{smallmatrix}
  0& 0 & -1\\
0 & -1 & 0\\
-1 & 0 & 0 \end{smallmatrix} \Bigg)w_{\ell}, 1\Bigg)
\Bigg( \begin{smallmatrix}
1 & \frac{B_1}{A_1} & \frac{C_2}{A_2} \\
0 & 1 & \frac{-B_2}{A_2} \\
0 & 0 & 1 \end{smallmatrix} \Bigg)
\Bigg( \begin{smallmatrix}
|4A_1| & 0 & 0 \\
0 & |\frac{A_2}{A_1}| & 0 \\
0 & 0 & \frac{1}{|4A_2|} \end{smallmatrix} \Bigg)^{-1}
\Bigg(\Bigg( \begin{smallmatrix}
\operatorname{sign}(A_1) & 0 & 0 \\
0 & \operatorname{sign}(\frac{A_2}{A_1}) & 0 \\
0 & 0 & \operatorname{sign}(\frac{1}{A_2}) \end{smallmatrix} \Bigg),\epsilon\Bigg)^{-1}.
\label{welliden}
\end{multline}
\item If $\gamma^\prime=\Bigg( \begin{smallmatrix}
0 & 0 & 1\\
1 & 0 & 0\\
0 & 1 & 0 \end{smallmatrix} \Bigg)
\Bigg( \begin{smallmatrix}
\frac{A_2}{B_1} & 0 & \frac{C_2}{4B_1}\\
0 & -4B_1 & \frac{4B_1B_2}{A_2}\\
0 & 0 & \frac{-1}{4A_2} \end{smallmatrix} \Bigg)$ and $\epsilon=(-B_1,-A_2)$, then
\begin{multline}(\gamma^\prime ,1)^{-1}(w_{\ell},1) =\\
 \Bigg(\Bigg( \begin{smallmatrix}
 0 & 1 & 0\\
0 & 0 & 1\\
1 & 0 & 0 \end{smallmatrix} \Bigg)
w_{\ell}, 1\Bigg)
\Bigg( \begin{smallmatrix}
1 & 0 & \frac{B_2}{A_2} \\
0 & 1 & \frac{C_2}{4A_2} \\
0 & 0 & 1 \end{smallmatrix} \Bigg)
\label{walpha12iden}\Bigg( \begin{smallmatrix}
|-4B_1| & 0 & 0 \\
0 & |\frac{A_2}{B_1}| & 0 \\
0 & 0 & |\frac{-1}{4A_2}| \end{smallmatrix} \Bigg)^{-1}
\Bigg(\Bigg( \begin{smallmatrix}
\operatorname{sign}(-B_1) & 0 & 0 \\
0 & \operatorname{sign}(\frac{A_2}{B_1}) & 0 \\
0 & 0 & \operatorname{sign}(\frac{-1}{A_2}) \end{smallmatrix} \Bigg),\epsilon\Bigg)^{-1}.
\end{multline}
\item If $\gamma^\prime=\Bigg( \begin{smallmatrix}
 0 & 1 & 0\\
0 & 0 & 1\\
1 & 0 & 0 \end{smallmatrix} \Bigg)
\Bigg( \begin{smallmatrix}
-4A_1 & -4B_1 & -C_1\\
 0& \frac{-1}{4B_2} & 0\\
 0& 0 & \frac{B_2}{A_1} \end{smallmatrix} \Bigg)$ and $\epsilon=-(A_1,-B_2)$, then
\begin{multline}(\gamma^\prime, 1)^{-1}(w_{\ell},1) =\\
\Bigg(\Bigg( \begin{smallmatrix}
 0 & 0 & 1\\
1 & 0 & 0\\
0 & 1 & 0 \end{smallmatrix} \Bigg)
w_{\ell}, 1\Bigg)
\Bigg( \begin{smallmatrix}
1 & \frac{C_1}{4A_1} & \frac{-B_1}{A_1} \\
0 & 1 & 0 \\
0 & 0 & 1 \end{smallmatrix} \Bigg)
\label{walpha21iden}\Bigg( \begin{smallmatrix}
|-4A_1| & 0 & 0 \\
0 & |\frac{B_2}{A_1}| & 0 \\
0 & 0 & |\frac{-1}{4B_2}| \end{smallmatrix} \Bigg)^{-1}
\Bigg(\Bigg( \begin{smallmatrix}
\operatorname{sign}(-A_1) & 0 & 0 \\
0 & \operatorname{sign}(\frac{B_2}{A_1}) & 0 \\
0 & 0 & \operatorname{sign}(\frac{-1}{B_2}) \end{smallmatrix} \Bigg),\epsilon\Bigg)^{-1}.
\end{multline}
\item If $\gamma^\prime=\Bigg( \begin{smallmatrix}
 0 & -1 & 0\\
-1 & 0 & 0\\
0  & 0 & -1 \end{smallmatrix} \Bigg)
\Bigg( \begin{smallmatrix}
\frac{4B_2}{C_1} & \frac{-C_2}{C_1} & 0\\
0 & \frac{1}{4B_2} & 0\\
0 & 0 & C_1 \end{smallmatrix} \Bigg)$ and $\epsilon=-(-C_1,-B_2)$, then
\begin{multline}(\gamma^\prime , 1)^{-1}(w_{\ell},1)=\\
\Bigg(\Bigg( \begin{smallmatrix}
 0 & -1 & 0\\
-1 & 0 & 0\\
0 & 0 & -1\end{smallmatrix} \Bigg)
w_{\ell}, 1\Bigg)
\Bigg( \begin{smallmatrix}
1 & 0 &0 \\
0 & 1 & \frac{-C_2}{4B_2} \\
0 & 0 & 1 \end{smallmatrix} \Bigg)
\label{walpha1iden}\Bigg( \begin{smallmatrix}
|C_1| & 0 & 0 \\
0 & |\frac{4B_2}{C_1}| & 0 \\
0 & 0 & |\frac{1}{4B_2}| \end{smallmatrix} \Bigg)^{-1}
\Bigg(\Bigg( \begin{smallmatrix}
\operatorname{sign}(C_1) & 0 & 0 \\
0 & \operatorname{sign}(\frac{B_2}{C_1}) & 0 \\
0 & 0 & \operatorname{sign}(\frac{1}{B_2}) \end{smallmatrix} \Bigg), \epsilon\Bigg)^{-1}.
\end{multline}
\item If $\gamma^\prime=\Bigg( \begin{smallmatrix}
 -1 & 0& 0\\
0& 0& -1\\
0& -1 & 0 \end{smallmatrix} \Bigg)
\Bigg( \begin{smallmatrix}
\frac{1}{C_2} & 0 & 0\\
0 & 4B_1 & C_1\\
0 & 0 & \frac{C_2}{4B_1} \end{smallmatrix} \Bigg)$ and $\epsilon=(B_1,-C_2)$, then
\begin{multline}(\gamma^\prime, 1)^{-1}(w_{\ell},1)=\\
\Bigg(\Bigg( \begin{smallmatrix}
-1 & 0 & 0\\
0& 0 & -1\\
0& -1 & 0 \end{smallmatrix} \Bigg)
w_{\ell}, 1\Bigg)
\Bigg( \begin{smallmatrix}
1 & \frac{C_1}{4B_1} & 0 \\
0 & 1 & 0\\
0 & 0 & 1 \end{smallmatrix} \Bigg)
\label{walpha2iden}\Bigg( \begin{smallmatrix}
|4B_1| & 0 & 0 \\
0 & |\frac{C_2}{4B_1}| & 0 \\
0 & 0 & |\frac{1}{C_2}| \end{smallmatrix} \Bigg)^{-1}
\Bigg(\Bigg( \begin{smallmatrix}
\operatorname{sign}(B_1) & 0 & 0 \\
0 & \operatorname{sign}(\frac{C_2}{B_1}) & 0 \\
0 & 0 & \operatorname{sign}(\frac{1}{C_2}) \end{smallmatrix} \Bigg),\epsilon\Bigg)^{-1}.
\end{multline}
\end{itemize}
\end{proposition}

\subsection{The Splitting}\label{splittingreview}

In \cite{M06}, Miller constructs a group homomorphism $S:\Gamma_{1}(4)\hookrightarrow\widetilde{\text{SL}}(3,\mathbb{R})$ such that $S(\gamma) = (\gamma,s(\gamma))$, where $s(\gamma)\in \{\pm1\}$. This map $S$ is a splitting of $\Gamma_{1}(4)$ into $\widetilde{SL}(3,\mathbb{R})$. Now we will describe a formula for $s$, which we may also call a splitting, in terms of Pl\"{u}cker coordinates. These results are proved in \cite{K18}.

\begin{theorem}[\cite{K18}] \label{theorem:PluckSplit} Let $\gamma\in\Gamma_{1}(4)$ with Pl\"{u}cker coordinates $(4A_1,4B_1,C_1,4A_2,4B_2,C_2)$ such that $A_1>0$, and $A_2/(A_1,A_2)\equiv 1\,(\text{mod } 2)$. Let $D=(A_1,A_2)$, $D_1 = (D,B_1)$, $D_2=D/D_1$, and let $\epsilon = \left(\frac{-1}{B_1/D_1}\right)$. Then
\begin{equation} \label{plucksplitting}
s(\gamma)=\left(\frac{\epsilon}{-A_1A_2}\right)\left(\frac{A_1/D}{A_2/D}\right)
\left(\frac{B_1/D_1}{A_1/D}\right)\left(\frac{4B_2/D_2}{\text{sign}(A_2)A_2/D}\right)\left(\frac{D_1}{C_1}\right)\left(\frac{D_2}{C_2}\right).\end{equation}

\end{theorem}

\begin{proposition}[\cite{K18}]\label{PluckSplitAll} Let $\gamma\in \Gamma_{1}(4)$ with Pl\"{u}cker Coordinates $(A_1,B_1,C_1,A_2,B_2,C_2)$. Then:

\begin{center}\begin{tabular}{|c|c|c|} \hline
\text{Cell} & $(A_1,B_1,C_1,A_2,B_2,C_2)$ & $s(\gamma)$\\ \hline
$B$ & $(0,0,-1,0,0,-1)$ & 1\\ \hline
$Bw_{\alpha_{1}}B$ & $(0,0,-1,0,B_2,C_2)$ & $\left(\frac{B_2}{-C_2}\right)$ \\ \hline
$Bw_{\alpha_{2}}B$ & $(0,B_1,C_1,0,0,-1)$ & $\left(\frac{-B_1}{-C_1}\right)$\\ \hline
$Bw_{\alpha_{1}}w_{\alpha_{2}}B$ & $(0,B_1,C_1,A_2,B_2,C_2)$ & $\left(\frac{A_2/B_1}{-C_2}\right)\left(\frac{-B_1}{-C_1}\right)$\\ \hline
$Bw_{\alpha_{2}}w_{\alpha_{1}}B$ & $(A_1,B_1,C_1,0,B_2,C_2)$ & $(-A_1,B_2)\left(\frac{-A_1/B_2}{-C_1}\right)\left(\frac{B_2}{-C_2}\right)$\\ \hline
$Bw_{\ell}B$ & $(A_1,B_1,C_1,A_2,B_2,C_2)$ & Theorem \ref{theorem:PluckSplit} \\ \hline
\end{tabular}
\end{center}
\end{proposition}

The splitting $s$ also satisfies identities induced from the maps described in Proposition \ref{prop:PluckerSym}.

\begin{proposition}[\cite{K18}] \label{CartSym}\label{doublecosetsym}\label{SignSym} Let $\gamma\in \Gamma_1(4)$ with Pl{\"u}cker coordinates $(4A_1,4B_1,C_1,4A_2,4B_2,C_2).$  Consider the involution $\psi:\gamma\mapsto w_{\ell}\,^{\intercal}\gamma^{-1}w_{\ell}^{-1}$. 
\begin{itemize}
\item If $A_1$ and $A_2$ are not equal to $0$, then $s(\psi(\gamma)) = (-A_1,-A_2)s(\gamma)$.
\item If $A_1,B_2\neq0$ and $A_2=0$, then $s(\psi(\gamma))=(-A_1,B_2)s(\gamma)$.
\item If $n\in \Gamma_{\infty}$, then $s(n\gamma)=s(\gamma n)=s(\gamma)$.
\item If $S_2=\mathrm{t}(1,-1,1)$ and $A_1,A_2\neq 0$, then $s(S_2\gamma S_2)=-\text{sign}(A_1A_2)s(\gamma)$.
\item If $S_3 =\mathrm{t}(1,1,-1)$ and $A_1,A_2\neq 0$, then $s(S_3\gamma S_3)=s(\gamma)$.
\end{itemize}
\end{proposition}

\textbf{Remark:} Using Proposition \ref{prop:PluckerSym}, Theorem \ref{theorem:PluckSplit}, Proposition \ref{PluckSplitAll}, and Proposition \ref{SignSym}, we may compute $s(\gamma)$ for any $\gamma\in\Gamma_{1}(4)$. Specifically, Proposition \ref{PluckSplitAll} provides a formula for any $\gamma\in\Gamma_{1}(4)$ that is not contained in the big Bruhat cell; Theorem \ref{theorem:PluckSplit} provides a formula for a subset of $\gamma\in\Gamma_{1}(4)$ in the big cell; Proposition \ref{prop:PluckerSym} and Proposition \ref{SignSym} allow the calculation of $s(\gamma)$ for any $\gamma\in\Gamma_{1}(4)$ in the big cell to be reduced to the subset covered in Theorem \ref{theorem:PluckSplit}.

We conclude this subsection by recalling the twisted multiplicativity of $s$.

\begin{proposition}[\cite{K18}] \label{twistmult} Let $A_1,\alpha_1\in\mathbb{Z}_{>0}$, $A_2,\alpha_2\in \mathbb{Z}$ such that $A_1,A_2$ are odd, \newline$(A_1A_2,\alpha_1\alpha_2)=1$, $A_1\alpha_1+A_2\alpha_2\equiv 0\text{ (mod }4)$, and $\frac{\alpha_2}{(\alpha_1,\alpha_2)}\equiv 1\text{ (mod }2)$. Let $\mu = \left(\frac{-1}{-A_1A_2}\right)$. Then with respect to the map from Proposition \ref{mult}, \newline$\varphi:S(A_1\alpha_1,A_2\alpha_2)\rightarrow S(A_1,\mu A_2)\times S(\alpha_1,-\mu\alpha_2)$, the following holds:
\begin{equation}s(\gamma) = s(\pi_1(\varphi(\gamma)))s(\pi_2(\varphi(\gamma)))\left(\frac{\alpha_2}{(\frac{-1}{A_1})A_1}\right)\left(\frac{\alpha_1}{A_2}\right),\end{equation}
where $\pi_i$ is the projection onto the $i$-th factor.
\end{proposition}


\subsection{Principal Series} \label{ssec:PSeries}

The following discussion establishes the preliminaries needed for the definition of the Eisenstein distribution. A complete treatment of automorphic distributions can be found in \cite{S00} and \cite{MS06}. Bate \cite{B} provides an explicit exposition of some of these ideas in the context of $\widetilde{\text{SL}}(2,\mathbb{R})$ principal series representations. We will solely be concerned with $\widetilde{\text{SL}}(3,\mathbb{R})$ principal series representations.


To construct a metaplectic principal series we employ the following representation. Let $(\phi,W)=(\phi,\mathbb{C}^{2})$ be the representation of $\widetilde{M}$ where $\phi:\widetilde{M}\rightarrow \text{SL}(2,\mathbb{C})$ is defined by
\begin{equation}\label{Q8rep}\phi\left(\left(\left(\begin{smallmatrix}
\epsilon_1 & 0 & 0\\
0 & \epsilon_2 & 0\\
0 & 0 & \epsilon_1\epsilon_2 \end{smallmatrix} \right),\pm1\right)\right)=\pm
\left(\begin{smallmatrix}
&-1\\
1&
\end{smallmatrix}\right)^{\frac{1-\epsilon_1}{2}}
\left(\begin{smallmatrix}
 -i&\\
&i
\end{smallmatrix}\right)^{\frac{1-\epsilon_2}{2}}.
\end{equation}

One can see that $\widetilde{M}$ is isomorphic to the quaternion group $Q_8=\{\pm1,\pm i,\pm j,\pm k\}$ and the representation $\phi$ is the unique irreducible two-dimensional complex representation of $Q_8$. The contragredient representation $(\phi^{\vee}, W^{\vee})$ can be realized as $(^{\intercal}\phi^{-1}, \mathbb{C}^{2})$. Under this identification, the natural pairing between $W$ and $W^{\vee}$ becomes the pairing $\langle\cdot,\cdot\rangle:\mathbb{C}^{2}\times \mathbb{C}^{2}\rightarrow \mathbb{C}$ defined by $\langle e_{i},e_{j}\rangle=\delta_{ij}$, where $e_{i}$ is the $i$-th standard basis vector in $\mathbb{C}^{2}$ and $\delta_{ij}$ is the Kronecker delta function.

Now we can define the principal series representations. Let $\lambda,\rho\in \frak{a}^\prime_{\mathbb{C}}$, where $\rho=(1,0,-1)$. Let
\begin{multline*}
\widetilde{V}_{\lambda,\phi}^{\infty}=\{f\in C^{\infty}(\widetilde{\text{SL}}(3,\mathbb{R}),W)|f(\tilde{g}\tilde{m}an_{-})=\text{exp}((\lambda-\rho)(H(a^{-1})))\phi(\tilde{m}^{-1})f(\tilde{g}),\\
\text{ for all } \tilde{m}an_{-}\in \widetilde{M}\widetilde{A}\widetilde{N}_{-}\},
\end{multline*}
\begin{multline*}
\widetilde{V}_{\lambda,\phi}=\{f\in L^2_{\text{loc}}(\widetilde{\text{SL}}(3,\mathbb{R}),W)|f(\tilde{g}\tilde{m}an_{-})=\text{exp}((\lambda-\rho)(H(a^{-1})))\phi(\tilde{m}^{-1})f(g),\\
\text{ for all } \tilde{m}an_{-}\in \widetilde{M}\widetilde{A}\widetilde{N}_{-}\},
\end{multline*}
\begin{multline*}
\widetilde{V}_{\lambda,\phi}^{-\infty}=\{f\in C^{-\infty}(\widetilde{\text{SL}}(3,\mathbb{R}),W)|f(\tilde{g}\tilde{m}an_{-})=\text{exp}((\lambda-\rho)(H(a^{-1})))\phi(\tilde{m}^{-1})f(g),\\
\text{ for all } \tilde{m}an_{-}\in \widetilde{M}\widetilde{A}\widetilde{N}_{-}\}.
\end{multline*}

The group $\widetilde{G}$ acts on each of the three preceding spaces by $\pi(h)(f)(g)=f(h^{-1}g)$. These spaces are called the smooth, locally $L^2$, and distributional principal series representation spaces, respectively. Note that $\widetilde{V}_{\lambda,\phi}^{\infty}\subset \widetilde{V}_{\lambda,\phi}\subset \widetilde{V}_{\lambda,\phi}^{-\infty}$. For $f_1=\left[\begin{array}{c}
f_{1,1}\\
f_{1,2}\end{array}\right]
\in \widetilde{V}_{-\lambda,^\intercal\phi^{-1}}$ and $f_2=\left[\begin{array}{c}
f_{2,1}\\
f_{2,2}\end{array}\right]
\in \widetilde{V}_{\lambda,\phi}$, define $(f_1\cdot f_2)(g)=f_{1,1}(g)f_{2,1}(g)+f_{1,2}(g)f_{2,2}(g)$.
We define the pairing $\langle\cdot,\cdot\rangle_{\lambda,\phi}:\widetilde{V}_{-\lambda,^\intercal\phi^{-1}}\times \widetilde{V}_{\lambda,\phi}\rightarrow \mathbb{C}$ by 
\begin{align} \label{metaCmptNoncmpt}
\langle f_1,f_2\rangle_{\lambda,\phi}=\int_{\widetilde{K}}(f_1\cdot f_2)(\tilde{k})d\tilde{k},
\end{align}
where $\widetilde{K}\cong \text{SU}(2)\cong \text{Spin}(3)$ and $dk$ is the Haar measure of $\widetilde{K}$. By a slight modification to \cite[Theorem 3]{LSL2}, if $h\in \widetilde{G}$, $f_1\in \widetilde{V}_{-\lambda,^\intercal\phi^{-1}}$, and $f_2\in \widetilde{V}_{\lambda,\phi}$, then $\langle\pi(h)f_1,\pi(h)f_2\rangle_{\lambda,\phi}=\langle f_1,f_2\rangle_{\lambda,\phi}$.

For us, distributions will be dual to smooth measures, and thus can be thought of as generalized functions in which the action of the distribution on the measure is given by integration of their product over the full space. Thus, the pairing can be extended to $\widetilde{V}_{\lambda,\phi}^{-\infty}$ on the right. Restriction from $\widetilde{V}_{-\lambda,^\intercal\phi^{-1}}$ to its smooth vectors results in a pairing $\widetilde{V}_{-\lambda,^\intercal\phi^{-1}}^{\infty}\times \widetilde{V}_{\lambda,\phi}^{-\infty}\rightarrow \mathbb{C}$. Under this pairing $\widetilde{V}_{\lambda,\phi}^{-\infty}$ may be identified with the dual of $\widetilde{V}_{-\lambda,^\intercal\phi^{-1}}^{\infty}$. This duality is to be understood in the context of topological vector spaces, thus some comments about topology are in order.

The map induced by restriction to $\widetilde{K}$ defines a vector space isomorphism between $\widetilde{V}_{-\lambda,^\intercal\phi^{-1}}^{\infty}$ and $C^{\infty}(\widetilde{K})$. The family of norms $||\partial^{\alpha}f||_u=\text{sup}_{k\in \widetilde{K}}\{|\partial^{\alpha}f(k)|\}$ defines a topology on $C^{\infty}(\widetilde{K})$ which can be transferred to $\widetilde{V}_{-\lambda,^\intercal\phi^{-1}}^{\infty}$ via the previous isomorphism. The dual $\widetilde{V}_{\lambda,\phi}^{-\infty}$ can be given the strong topology \cite[\S 19]{T67}. With respect to these topologies $\widetilde{V}_{\lambda,\phi}^{-\infty}$ can be identified with the continuous dual of $\widetilde{V}_{-\lambda,^\intercal\phi^{-1}}^{\infty}$. Additionally, $\widetilde{V}_{\lambda,\phi}^{\infty}$ is dense in $\widetilde{V}_{\lambda,\phi}^{-\infty}$, and sequential convergence in $\widetilde{V}_{\lambda,\phi}^{-\infty}$ with respect to the strong topology is equivalent to sequential convergence with respect to the weak topology \cite[\S 34.4]{T67}.

The pairing just described focuses on the compact model of the principal series representations. The Eisenstein distribution considered in this paper will be more amenable to study using the non-compact model of the principal series representation which we describe presently.

Let $w\in W$. As $w\widetilde{N}\widetilde{B}_{-}$ is open and dense in $\widetilde{G}$, restriction from $\widetilde{G}$ to $w\widetilde{N}$ defines an injection $\widetilde{V}_{\lambda,\phi}^{\infty}\hookrightarrow C^{\infty}(w\widetilde{N})$ and the pairing is compatible with this injection in the following sense. 
Let $F:\widetilde{G}\rightarrow \mathbb{C}$, be a smooth function such that $F(gb_{-}) = e^{2\rho (H(b_{-}))}F(g)$. Then, by a slight modification of Consequence 7 in \cite{K86}, 
\begin{equation}\label{KNintegral} \int_{\widetilde{K}}F(k)dk = \int_{\widetilde{N}}F(wn)dn.
\end{equation}
For $b_{-}\in \widetilde{B}_{-}$ we have 
$(f_1\cdot f_2)(gb_{-})=e^{2\rho(H(b_{-}))}(f_1\cdot f_2)(g)$.
This identity allows us to apply Equation (\ref{KNintegral}) to establish a bridge between the pairings of principal series in the compact and noncompact pictures. Specifically, 
\begin{equation}\langle f_1,f_2\rangle_{\lambda,\phi}=\int_{\widetilde{N}}(f_1\cdot f_2)((w,1)n)dn,\label{metaCmptNCmpt}
\end{equation}
and so the pairing can be realized as an integration over the non-compact space $\widetilde{N}$.

The element $\tilde{\tau}\in \widetilde{V}_{\lambda,\phi}^{-\infty}$, characterized by
\begin{equation}\label{metatau}\tilde{\tau}\left((w_{\ell},1)
\left( \begin{smallmatrix}
1 & x & z\\
0 & 1 & y\\
0 & 0 & 1 \end{smallmatrix} \right)\tilde{m}an_{-}\right)
=\text{exp}((\lambda-\rho)(H(a^{-1})))\phi(\tilde{m}^{-1})\left[\begin{array}{c}
\delta_{(0,0,0)}(x,y,z)\\
0\end{array}\right],\end{equation}
will be used to construct a metaplectic Eisenstein distribution on $\widetilde{\text{SL}}(3,\mathbb{R})$.

\begin{proposition} \label{metatauprops} Let $\tilde{\tau}\in \widetilde{V}_{\lambda,\phi}^{-\infty}$ be as above. Then:
\begin{enumerate}
\item $\tilde{\tau}$ is right $\widetilde{N}_{-}$-invariant.
\item $\text{supp}(\tilde{\tau})= (w_{\ell},1)\widetilde{B}_{-}=\widetilde{B}(w_{\ell},1)$.
\item $\tilde{\tau}$ is left $\widetilde{N}$-invariant.
\end{enumerate}
\end{proposition}
\textbf{Proof:} The first two properties follow immediately from the definition of $\tilde{\tau}$. For the final claim let $f=\left[\begin{smallmatrix}
f_1\\
f_2
\end{smallmatrix}\right]\in \widetilde{V}_{-\lambda,^\intercal\phi^{-1}}^{\infty}$. By the definition of $\tilde{\tau}$, we have $\langle f,\tilde{\tau}\rangle_{\lambda,\phi}=f_1(w_{\ell})$.
 On the other hand, $\langle f,\pi(n)\tilde{\tau}\rangle_{\lambda,\phi} = \langle \pi(n^{-1})f,\tilde{\tau}\rangle_{\lambda,\phi} = f_1(nw_{\ell}) = f_1(w_{\ell})$,
 as $n(w_{\ell},1)=(w_{\ell},1)n_{-}$, where $n_{-}\in \widetilde{N}_{-}=(w_{\ell},1)^{-1}\widetilde{N}(w_{\ell},1)$. Thus, $\pi(n)\tilde{\tau}=\tilde{\tau}$.\EndProof

Now we can define the metaplectic Eisenstein distribution as 
\begin{equation}\label{metaEisen}\tilde{E}(\tilde{g},\lambda)
 = \sum_{\gamma\in \Gamma_{\infty}\backslash\Gamma}\pi(S(\gamma)^{-1})\tilde{\tau}(\tilde{g})\in \widetilde{V}_{\lambda,\phi}^{-\infty},\end{equation}
where 
$\tilde{\tau}$ is as in line (\ref{metatau}) and $\lambda=(\lambda_1,\lambda_2,\lambda_3)\in \mathbb{C}^3$ such that $\lambda_1+\lambda_2+\lambda_3=0$. Note that $\tilde{E}$ is well defined since $\tilde{\tau}$ is left $\widetilde{N}$-invariant. Furthermore, standard arguments show that the sum is convergent when the real parts of $\lambda_1-\lambda_2$ and $\lambda_2-\lambda_3$ are sufficiently large. 



\section{Exponential Sums}\label{sec:ExpSums}

\subsection{Preliminaries} \label{sec:ExpSumPrelim}

Let $A_1,A_2\in\mathbb{Z}_{\neq0}$. This section begins a study of the exponential sums, $\Sigma(A_1,A_2;m_1,m_2)$ (defined below), that appear as the coefficients of the Dirichlet series that make up the Fourier coefficients of the metaplectic Eisenstein distribution associated with the big cell. The primary focus of this section is to reduce the general computation of $\Sigma(A_1,A_2;m_1,m_2)$ to the case where $A_1$ and $A_2$ are prime powers.

We begin by describing how the symmetries of $s$, described in subsection \ref{splittingreview}, affect
\begin{equation}\label{SigmaA1A2}\Sigma(A_1,A_2;m_1,m_2) \stackrel{\text{def}}{=} \sum_{\gamma\in \Gamma_{\infty}\backslash\mathbb{S}(A_1,A_2)/\Gamma_{\infty}}s(\gamma)e^{2\pi i(m_1\frac{B_1}{A_1}+m_2\frac{B_2}{A_2})}.\end{equation}

\begin{proposition} \label{ExpSym} Let $A_1,A_2\in\mathbb{Z}_{\neq0}$ and let $m_1,m_2\in\mathbb{Z}$. Then:
\begin{enumerate}
\item $\Sigma(A_1,A_2;m_1,m_2)=\Sigma(-A_1,-A_2;m_1,-m_2)$.
\item $\Sigma(A_1,A_2;m_1,m_2)=(-A_1,-A_2)\Sigma(A_2,A_1;-m_2,-m_1)$.
\item If $A_1A_2>0$, then $\Sigma(A_1,A_2;0,0)=0$.
\end{enumerate}
\end{proposition}

\textbf{Proof:} This follows from the definition of $\Sigma(A_1,A_2;m_1,m_2)$ and propositions \ref{prop:PluckerSym} and \ref{CartSym}.\EndProof

We can also apply the twisted multiplicativity of Proposition \ref{twistmult} to study $\Sigma(A_1,A_2;m_1,m_2)$.

\begin{proposition} \label{exptwistmult} Let $A_1,\alpha_1\in\mathbb{Z}_{>0}$, $A_2,\alpha_2\in \mathbb{Z}$ such that $A_1,A_2$ are odd, $(A_1A_2,\alpha_1\alpha_2)=1$, and $\alpha_2$ is divisible by fewer powers of 2 than $\alpha_1$. Let $\mu = \left(\frac{-1}{-A_1A_2}\right)$. Then with respect to the map from Proposition \ref{mult}, $S(A_1\alpha_1,A_2\alpha_2)\stackrel{\varphi}{\rightarrow} S(A_1,\mu A_2)\times S(\alpha_1,-\mu\alpha_2)$, the following holds:
\begin{multline*}
\Sigma(A_1\alpha_1,A_2\alpha_2;m_1,m_2) \\= \left(\frac{\alpha_2}{(\frac{-1}{A_1})A_1}\right)\left(\frac{\alpha_1}{A_2}\right)\Sigma(A_1,\mu A_2;(\alpha_1)^{-1}_{\text{mod }A_1}m_1,(\alpha_2)^{-1}_{\text{mod }A_2}m_2)\\
\times\Sigma(\alpha_1,-\mu\alpha_2;(A_1)^{-1}_{\text{mod }\alpha_1}m_1,\left(\frac{-1}{A_1}\right)(A_2)^{-1}_{\text{mod }\alpha_2}m_2).
\end{multline*}
\end{proposition}

\textbf{Proof:} This follows from Proposition \ref{twistmult}. \EndProof

The function $\Sigma(A_1,A_2;m_1,m_2)$ also exhibits a multiplicativity in the variables $m_1,m_2$.

\begin{proposition} \label{twistmultindex} Let $A_1,A_2,m_1,m_2,c_1,c_2\in\mathbb{Z}$ such that, $A_1>0$, $A_2,m_1,m_2\neq 0$, and $(c_1c_2,A_1A_2)=1$. Then $\Sigma(A_1,A_2;c_1m_1,c_2m_2) = \left(\frac{c_1}{A_1}\right)\left(\frac{c_2}{A_2}\right)\Sigma(A_1,A_2;m_1,m_2)$.
\end{proposition}

\noindent Proposition \ref{twistmultindex} is the byproduct of Proposition \ref{exptwistmult} and our computations in subsections \ref{ssec:Sigma2} and \ref{ssec:Sigma4}, so no direct proof will be given.

Proposition \ref{exptwistmult} and Proposition \ref{twistmultindex} show that the exponential sums $\Sigma(A_1,A_2;m_1,m_2)$ are built from those of the form $\Sigma(\pm p^k,\pm p^\ell;p^{r_1},p^{r_2})$, where $p$ is a prime.


\subsection{Explicit Description of Double Cosets: $A_1$, $A_2$ Odd} \label{sec:DesofDoubleCosets}

The sum $\Sigma(A_1,A_2;m_1,m_2)$ is indexed by the set $\Gamma_{\infty}\backslash\mathbb{S}(A_1,A_2)/\Gamma_{\infty}\cong S(A_1,A_2)$. In this section we provide a description of the sets $S(A_1,A_2)$ when $A_1,A_2\neq0$. We begin with a few simple observations.

Note that as $C_j\equiv -1$ (mod $4$), it follows that $A_1\equiv -A_2$ (mod $4)$. Thus, if $A_1\not\equiv -A_2$ (mod $4)$, then $S(A_1,A_2) = \emptyset$; thus $\Sigma(A_1,A_2;m_1,m_2)=0$. By Proposition \ref{ExpSym}, it suffices to consider the case in which $A_1>0$; by Proposition \ref{mult} it suffices to study $S(p^k,\pm p^l)$, where $p$ is prime. However, a description of $S(A_1,A_2)$, with $(A_1,A_2)=1$, will be included as it provides a clean presentation of the boundary case $A_1=p^k$, $A_2=\pm1$. The proofs of the results of this subsection are straightforward and will be omitted. 

\begin{proposition}\label{relprime}
Let $A_1,A_2\in\mathbb{Z}$ such that $(A_1,A_2) = 1$, $A_1+A_2\equiv 0$ (mod $4$). Then $S(A_1,A_2)$ consists precisely of the elements:
\[(A_1,B_1,\frac{A_1C_2+4B_1B_2}{-A_2},A_2,B_2,C_2),\]
where: \begin{itemize}
\item $0\leq B_1<A_1$ and $(B_1,A_1)=1$;
\item $0\leq -B_2 <|A_2|$ and $(B_2,A_2)=1$;
\item $C_2\equiv -(A_1)^{-1}(4B_1B_2)$ (mod $A_2$), $C_2\equiv -1$ (mod $4$), and $0\leq \text{sign}(A_2)C_2<4|A_2|$.
\end{itemize}
In particular, $|S(A_1,A_2)|=\phi(A_1)\phi(A_2)$.
\end{proposition}


The description of $S(p^k,\pm p^\ell)$ will be broken into several cases: $k>\ell>0$, $k=\ell>0$, and $0<k<\ell$. However, Proposition \ref{prop:PluckerSym} reveals that the last case is redundant. 

\begin{proposition} \label{+-p,k>l>0}
Let $p$ be an odd prime, $\mu=\pm1$, and let $k>\ell>0$ be integers such that $p^k\equiv -\mu p^\ell\,(\text{mod }4)$. Then $S(p^k,\mu p^{\ell})$ consists precisely of the elements:
\begin{enumerate}
\item$(p^k,B_1,p^{k-\ell}C_2,\mu p^{\ell},0,C_2)$
where:
\begin{itemize}
\item $0\leq B_1<p^k$ and $(B_1,p^k)=1$;
\item $0\leq \mu C_2<4p^l$, $C_2\equiv -1 \, (\text{mod }4)$, and $(C_2,p^{\ell})=1$.
\end{itemize}
\item $(p^k,b_1p^{\ell},p^kC_2+4b_1b_2,\mu p^{\ell},b_2,C_2)$ where:
\begin{itemize}
\item $0<b_1<p^{k-\ell}$, $(b_1,p^k)=1$;
\item $0<\mu b_2<p^{\ell}$, $(b_2,p^{k-\ell})=1$;
\item $0\leq \mu C_2<4p^l$, and $C_2\equiv -1 \, (\text{mod }4)$.
\end{itemize}
\item $(p^k,b_1p^i,p^kC_2+4b_1b_2,\mu p^{\ell},b_2p^{\ell-i},C_2)$ where:
\begin{itemize}
\item $0<i<l$;
\item $0<b_1<p^{k-i}$, $(b_1,p^k)=1$;
\item $0<\mu b_2<p^{k-\ell+i}$, $(b_2,p^{k-\ell+i})=1$;
\item $0\leq \mu C_2<4p^l$, $C_2\equiv -1 \, (\text{mod }4)$, and $(C_2,p^{\ell})=1$.
\end{itemize}
\end{enumerate}
\end{proposition}

\begin{proposition} \label{+-p,k>0}
Let $p$ be an odd prime and let $k\in\mathbb{Z}_{>0}$. Then $S(p^k,-p^k)$ consists precisely of the elements:
\begin{enumerate}
\item$(p^k,0,C_2,-p^k,0,C_2)$
where:
\begin{itemize}
\item $0\leq -C_2<4p^k$, $C_2\equiv -1 \, (\text{mod }4)$, and $(C_2,p^k)=1$.
\end{itemize}
\item$(p^k,0,C_2,-p^k,B_2,C_2)$
where:
\begin{itemize}
\item $0< -B_2<p^k$;
\item $0\leq -C_2<4p^k$, $C_2\equiv -1 \, (\text{mod }4)$, and $(C_2,p^k)=1$.
\end{itemize}
\item$(p^k,B_1,C_2,-p^k,0,C_2)$ where:
\begin{itemize}
\item $0<B_1<p^k$;
\item $0\leq -C_2<4p^k$, $C_2\equiv -1 \, (\text{mod }4)$, and $(C_2,p^k)=1$.
\end{itemize}
\item $(p^k,b_1p^i,C_2+4b_1b_2,-p^k,b_2p^j,C_2)$ where:
\begin{itemize}
\item $0<i<k$, $0<j<k$, and $k=i+j$;
\item $0<b_1<p^{k-i}$, and $(b_1,p^k)=1$;
\item $0<-b_2<p^{k-j}$, and $(b_2,p^k)=1$;
\item $0\leq -C_2<4p^k$, $C_2\equiv -1 \, (\text{mod }4)$, and $(C_2,p^k)=1$;
\item $C_2+4b_1b_2\not\equiv\,0$ (mod $p$).
\end{itemize}
\item $(p^k,b_1p^i,C_2+4b_1b_2,-p^k,b_2p^j,C_2)$ where:
\begin{itemize}
\item $0<i<k$, $0<j<k$, and $k<i+j$;
\item $0<b_1<p^{k-i}$, and $(b_1,p^k)=1$;
\item $0<-b_2<p^{k-j}$, and $(b_2,p^k)=1$;
\item $0\leq -C_2<4p^l$, $C_2\equiv -1 \, (\text{mod }4)$, and $(C_2,p^k)=1$.
\end{itemize}

\end{enumerate}
\end{proposition}

\subsection{Explicit Description of Double Cosets: $A_1$, $A_2$ Even} \label{sec:DesofDoubleCosetsEven}

It remains to consider $S(A_1,A_2)$ where $A_1=2^k$ and $A_2=\pm2^\ell$. Before describing $S(2^k,\pm2^\ell)$ a few comments are in order. First, as $A_i$ is a power of $2$ and $C_i$ must be odd, the condition $(A_i,B_i,C_i)=1$  is vacuous. Second, when $k>\ell$ the equation $2^{k-\ell}C_2+2^{2-\ell}B_1B_2\pm C_1=0$,
implies that $\ell\geq2$. Once again the following results are straightforward and the proofs are omitted.


\begin{proposition}\label{+pm2,k>l>0} Let $k>\ell\geq 2$ be integers and let $\mu=\pm1$. Then $S(2^k,\mu 2^{\ell})$ consists precisely of the elements:
\begin{enumerate}
\item$(2^k,2^ib_1,2^{k-\ell}C_2+b_1b_2,\mu2^{\ell},2^jb_2,C_2)$
where:
\begin{itemize}
\item $0\leq i <k$, and $0\leq j <\ell $ such that $i+j=\ell-2$;
\item $0\leq \mu C_2<4(2^{\ell})$, $C_2\equiv -1\, (\text{mod }4)$;
\item $0\leq b_1<2^{k-i}$ such that $(2,b_1)=1$;
\item $0\leq \mu b_2<2^{\ell-j}$ such that $(2,b_2)=1$, and $b_1b_2\equiv \mu+2^{k-\ell}\,(\text{mod }4)$.
\end{itemize}
\end{enumerate}
\end{proposition}

\begin{proposition}\label{+-2,k>0} Let $k\in\mathbb{Z}_{>0}$. Then $S(2^k,-2^{k})$ consists precisely of the elements:
\begin{enumerate}
\item$(2^k,0,C_2,-2^k,0,C_2)$ where:
\begin{itemize}
\item $0\leq -C_2<4(2^k)$, $C_2\equiv -1\, (\text{mod }4)$.
\end{itemize}
\item$(2^k,0,C_2,-2^k,B_2,C_2)$ where:
\begin{itemize}
\item $0\leq -C_2<4(2^k)$, $C_2\equiv -1\, (\text{mod }4)$;
\item $0<-B_2<2^{k}$.
\end{itemize}
\item$(2^k,B_1,C_2,-2^k,0,C_2)$ where:
\begin{itemize}
\item $0\leq -C_2<4(2^k)$, $C_2\equiv -1\, (\text{mod }4)$;
\item $0<B_1<2^{k}$.
\end{itemize}
\item$(2^k,2^ib_1,C_2+2^{2+i+j-k}b_1b_2,-2^k,2^jb_2,C_2)$
where:
\begin{itemize}
\item $0\leq i <k$, and $0\leq j <k$ such that $i+j\geq k$;
\item $0\leq -C_2<4(2^k)$, $C_2\equiv -1\, (\text{mod }4)$;
\item $0\leq b_1<2^{k-i}$ such that $(2,b_1)=1$;
\item $0\leq -b_2<2^{k-j}$ such that $(2,b_2)=1$.
\end{itemize}
\end{enumerate}
\end{proposition}

\begin{proposition}\label{++2,k>0} Let $k\in\mathbb{Z}_{>0}$. Then $S(2^k,2^{k})$ consists precisely of the elements:
\begin{enumerate}
\item$(2^k,2^ib_1,-C_2-2b_1b_2,2^k,2^jb_2,C_2)$
where:
\begin{itemize}
\item $0\leq i <k$, and $0\leq j <k$ such that $i+j =k-1$;
\item $0\leq C_2<4(2^k)$, $C_2\equiv -1\, (\text{mod }4)$;
\item $0\leq b_1<2^{k-i}$ such that $(2,b_1)=1$;
\item $0\leq b_2<2^{k-j}$ such that $(2,b_2)=1$.
\end{itemize}
\end{enumerate}
\end{proposition}

\subsection{Exponential Sum: $\Sigma(p^k,\pm p^l;m_1,m_2)$} \label{ssec:Sigma2}

In this subsection and the next we use the description of the set $S(A_1,A_2)$ contained in sections \ref{sec:DesofDoubleCosets} and \ref{sec:DesofDoubleCosetsEven}, and the formula for the splitting (Theorem \ref{theorem:PluckSplit}) to evaluate the exponential sum $\Sigma(A_1,A_2;m_1,m_2)$ defined on line (\ref{SigmaA1A2}). In particular, by Proposition \ref{exptwistmult} the computations can be reduced to the case where $A_1$ and $A_2$ are powers of a fixed prime. We begin with $p$ an odd prime.


The proofs of the next two propositions, \ref{sigma+-p k,0,m1,m2} and \ref{sigma+-p k>l>0,m1,m2}, are similar to the proof of Proposition \ref{sigma+-p,k>0,m1,m2} and will be omitted.
\begin{proposition} \label{sigma+-p k,0,m1,m2}Let $p$ be an odd prime, let $\mu=\pm1$, and let $k,m_1,m_2\in \mathbb{Z}$ such that $k>0$. 
Then
\begin{align*}\Sigma(1,\mu;m_1,m_2)=&\,\delta_{\mu=-1},\\
\text{and }\Sigma(p^k,\mu;m_1,m_2) =&
\begin{cases}
0,\hspace{2cm}\text{ if } p^k\equiv \mu  (\text{mod }4);\\
g(p^k,m_1,p^k),\text{ if } p^k\equiv -\mu (\text{mod }4).
\end{cases}\end{align*}

\end{proposition}

\begin{proposition} \label{sigma+-p k>l>0,m1,m2}Let $p$ be an odd prime, let $\mu=\pm1$, and let $k,\ell,m_1,m_2\in \mathbb{Z}$ such that $k>\ell>0$.
\begin{itemize}
\item If $p^{k-\ell}\equiv \phantom{-}\mu$ (mod $4$), then $\Sigma(p^k,\mu p^\ell;m_1,m_2)=0$.
\item If $p^{k-\ell}\equiv -\mu$ (mod $4$), then $\Sigma(p^k,\mu p^\ell;m_1,m_2)=$
\begin{equation*}p^{\ell}g(p^\ell, m_2,p^\ell)g(p^k,m_1,p^{k-\ell})
+\phi(p^l)\sum_{
\begin{smallmatrix}
0\leq i<l,\\
 i\equiv l\,(\text{mod }2)
 \end{smallmatrix}}
 g(p^i, m_2,p^i)g(p^k,m_1,p^{k-i}).
 \end{equation*}
 \end{itemize}
\end{proposition}
\begin{proposition} \label{sigma+-p,k>0,m1,m2}Let $p$ be an odd prime and let $k,m_1,m_2\in \mathbb{Z}$ such that $k>0.$  Then

\begin{align*}\Sigma(p^k,-p^k;m_1,m_2) =\sum_{i=1}^{k-1}p^{k-i}g(p^i,m_2,p^i)&g(p^{k-i},-m_2,p^i)g(p^k,m_1,p^{k-i})\\
+\delta_{2|k}&\phi(p^k)(\delta_{p^{k}|m_1}p^k+\delta_{p^{k}|m_2}p^k-\delta_{p^{k-1}|m_1}p^{k-1}).
\end{align*}

\end{proposition}
\noindent\textbf{Proof:} The summands of $\Sigma(p^k,-p^k;m_1,m_2)$ will be grouped into five cases in accordance with Proposition \ref{+-p,k>0}. However, it will be more convenient to include case 1 in both case 2 and case 3, and then subtract those terms that are double counted.

Case 2:  $B_1=0.$ For $\gamma\in S(p^k,-p^k)$ such that $B_1=0$, equation (\ref{plucksplitting}) implies that $s(\gamma) = (\frac{p^k}{C_1})=(\frac{p^k}{-C_1})=(\frac{-C_1}{p^k})=(\frac{-C_2}{p^k})$. 
By Proposition \ref{+-p,k>0} the sum under consideration becomes
\begin{equation}\label{kkoddeq2}\delta_{2|k}\delta_{p^k|m_2}p^k\phi(p^k).\end{equation}


Case 3:  $B_2=0.$ For $\gamma\in S(p^k,-p^k)$ such that $B_2=0$, equation (\ref{plucksplitting}) implies that $s(\gamma) = (\frac{D_1}{C_1})(\frac{D_2}{C_2})=(\frac{D}{C_2})=(\frac{p^k}{C_2})=(\frac{-C_2}{p^k})$, where the second equality follows from the identity $A_1C_2+4B_1B_2+C_1A_2=0$. By Proposition \ref{+-p,k>0}, the sum corresponding to this case is given by
\begin{equation}\label{kkoddeq3}\delta_{2|k}\delta_{p^k|m_1}p^k\phi(p^k).\end{equation}


Remember that we included case 1, $B_1=B_2=0$, in both case 2 and case 3. To compensate for this we subtract $\delta_{2|k}\phi(p^k)$.

Case 4: $B_1B_2\neq 0$, $i+j=k$,  where $B_1 =p^ib_1$ and $B_2 = p^jb_2$ such that $(b_1b_2,p) = 1.$  As $0\leq i,j<k$ it follows that $0<i,j$.  We apply equation (\ref{plucksplitting}) to see that if $\gamma\in S(p^k,-p^k)$ and $\gamma$ satisfies the conditions of case 4, then $s(\gamma)=(\frac{p^i}{C_1})(\frac{p^{k-i}}{C_2}) =(\frac{-C_1}{p^i})(\frac{-C_2}{p^{k-i}})
=(\frac{-C_2-4b_1b_2}{p^i})(\frac{-C_2}{p^{k-i}})$.
Now we can sum over the elements of this case described in Proposition \ref{+-p,k>0} with $i$ fixed to get
\begin{equation}\label{expsum1}\sum_{\begin{smallmatrix}
C_2\in (\mathbb{Z}/4p^k\mathbb{Z})^\times\\
C_2\equiv -1 \text{(mod }4)
\end{smallmatrix}}\sum_{b_1\in(\mathbb{Z}/p^{k-i}\mathbb{Z})^\times}\sum_{b_2\in(\mathbb{Z}/p^{i}\mathbb{Z})^\times}\left(\frac{-C_2-4b_1b_2}{p^i}\right)\left(\frac{-C_2}{p^{k-i}}\right)e^{2\pi i(m_1\frac{b_1}{p^{k-i}}-m_2\frac{b_2}{p^i})}.
\end{equation}
To simplify this expression consider the sum indexed by $b_2$ in line (\ref{expsum1}). By Lemma \ref{lemma:gauss},
\begin{multline}\label{aaabbb}
\sum_{b_2}\left(\frac{-C_2+(-4b_1)(b_2)}{p}\right)^ie^{2\pi i(-m_2\frac{b_2}{p^i})}\\
=e^{2\pi i(\frac{-m_2\overline{4b_1}(-C_2)}{p^i})}\left(\frac{b_1}{p}\right)^ig(p^i,m_2,p^i)-\left(\frac{-C_2}{p}\right)^i\delta_{p^{i-1}|m_2}p^{i-1}.\end{multline}
Next we consider the sum over $b_2$ and $C_2$ in line (\ref{expsum1}) and apply equation (\ref{aaabbb}) to get
\begin{multline}
\left(\frac{b_1}{p}\right)^i g(p^i,m_2,p^i)p^{k-i}\sum_{\begin{smallmatrix}
C_2\in (\mathbb{Z}/4p^i\mathbb{Z})^\times\\
C_2\equiv -1 \text{(mod }4)
\end{smallmatrix}}\left[\left(\frac{-C_2}{p^{k-i}}\right)e^{2\pi i(\frac{-m_2\overline{4b_1}(-C_2)}{p^i})}-\left(\frac{-C_2}{p}\right)^k\delta_{p^{i-1}|m_2}p^{i-1}\right]\\
=\left(\frac{b_1}{p}\right)^i g(p^i,m_2,p^i)
p^{k-i}\sum_{c\in (\mathbb{Z}/p^{i}\mathbb{Z})^\times}\left(\frac{c}{p^{k-i}}\right)e^{2\pi i(\frac{-m_2\overline{4b_1}c}{p^i})}
-\phi(p^k)\delta_{2|k}\delta_{p^{i-1}|m_2}p^{i-1}, \label{expsum+-kkeqn1}
\end{multline}
where the previous equality follows from an application of the Chinese Remainder Theorem. We simplify the result to get
\begin{align}\nonumber
(\ref{expsum+-kkeqn1})=&\left(\frac{b_1}{p}\right)^k g(p^i,m_2,p^i)
p^{k-i}\sum_{c\in (\mathbb{Z}/4p^{i}\mathbb{Z})^\times}\left(\frac{c}{p^{k-i}}\right)e^{2\pi i(\frac{-m_2c}{p^i})}
-\phi(p^k)\delta_{2|k}\delta_{p^{i-1}|m_2}p^{i-1}\\
=&\left(\frac{b_1}{p}\right)^k g(p^i,m_2,p^i)
p^{k-i}g(p^{k-i},-m_2,p^{i})
-\phi(p^k)\delta_{2|k}\delta_{p^{i-1}|m_2}p^{i-1}.\label{aabbb}\end{align}
Now consider (\ref{expsum1}) in its entirety and apply the simplifications that resulted in line (\ref{aabbb}) to see that (\ref{expsum1}) is equal to 
\begin{align}
\nonumber&\sum_{b_1}\left(\frac{b_1}{p}\right)^kg(p^i,m_2,p^i)p^{k-i}g(p^{k-i},-m_2,p^i)e^{2\pi i(m_1\frac{b_1}{p^{k-i}})}
-\phi(p^k)\delta_{2|k}\delta_{p^{i-1}|m_2}p^{i-1}e^{2\pi i(m_1\frac{b_1}{p^{k-i}})}\\
=&g(p^i,m_2,p^i)p^{k-i}g(p^{k-i},-m_2,p^i)g(p^k,m_1,p^{k-i})
-g(p^2,m_1,p^{k-i})\delta_{2|k}\phi(p^k)\delta_{p^{i-1}|m_2}p^{i-1}.\label{kkoddeq4}\end{align}
This must be summed over $0<i<k$.

Case 5: $B_1B_2\neq 0$, $i+j>k$,  where $B_1 =p^ib_1$ and $B_2 = p^jb_2$ such that $(b_1b_2,p) = 1.$ 
Using equation (\ref{plucksplitting}), it follows that if $\gamma\in S(p^k,-p^k)$ satisfies the conditions of case 5, then $s(\gamma) = (\frac{p^i}{C_1})(\frac{p^{k-i}}{C_2})=(\frac{p^i}{-C_1})(\frac{p^{k-i}}{-C_2})
=(\frac{-C_1}{p^i})(\frac{-C_2}{p^{k-i}})=(\frac{-C_2}{p^{k}})$.
The last equality follows as $C_1\equiv C_2$ (mod $p$). Thus the sum under consideration becomes
\begin{multline*}\sum_{b_1\in (\mathbb{Z}/p^{k-i}\mathbb{Z})^\times}\sum_{b_2\in (\mathbb{Z}/p^{k-j}\mathbb{Z})^\times}\sum_{\begin{smallmatrix}
-C_2\in (\mathbb{Z}/4p^k\mathbb{Z})^\times\\
-C_2\equiv 1 \text{(mod }4)
\end{smallmatrix}}(\frac{-C_2}{p^k})e^{2\pi i(m_1\frac{b_1}{p^{k-i}}-m_2\frac{b_2}{p^{k-j}})}\\
=\delta_{2|k}\phi(p^k)g(p^2,m_1,p^{k-i})g(p^2,m_2,p^{k-j}).
\end{multline*}
This must be summed over $0<i<k,$ and $k-i+1\leq j\leq k-1$, and then can be simplified using equation (\ref{Gauss5}). Specifically, 

\begin{align}
\nonumber\delta_{2|k}\sum_{0<i<k}\sum_{1\leq k- j\leq i-1}&\phi(p^k)g(p^2,m_1,p^{k-i})g(p^2,m_2,p^{k-j})\\
\nonumber=&\delta_{2|k}\phi(p^k)\sum_{0<i<k}g(p^2,m_1,p^{k-i})(-1+\delta_{p^{i-1}|m_2}p^{i-1})\\
\nonumber=&\delta_{2|k}\phi(p^k)(-\sum_{0<i<k}g(p^2,m_1,p^{k-i})+\sum_{0<i<k}g(p^2,m_1,p^{k-i})\delta_{p^{i-1}|m_2}p^{i-1})\\
\label{kkoddeq5}=&\delta_{2|k}\phi(p^k)(1-\delta_{p^{k-1}|m_1}p^{k-1}+\sum_{0<i<k}g(p^2,m_1,p^{k-i})\delta_{p^{i-1}|m_2}p^{i-1}).
\end{align}

The final result follows once we add (\ref{kkoddeq2}), (\ref{kkoddeq3}), (\ref{kkoddeq4}), and (\ref{kkoddeq5}) together, remembering to sum (\ref{kkoddeq4}) over $0<i<k$ and to remove what was double counted in cases 1 and 2.  \EndProof


\subsection{Exponential Sum: $\Sigma(2^k,\pm2^l;m_1,m_2)$}\label{ssec:Sigma4}

In this subsection we consider the bad prime $p=2$.  Recall that an explicit parameterization of the index of summation of \newline$\Sigma(2^k,\pm2^l;m_1,m_2)$ is described in Section \ref{sec:DesofDoubleCosetsEven}.

\begin{proposition} \label{sigma 2,k>l>0,m1,m2} Let $k,\ell\in \mathbb{Z}$ such that $k>\ell\geq0$, let $m_1,m_2\in\mathbb{Z}$, and let $\mu=\pm1$.
\begin{itemize}
\item If $\ell=0,1$, then $\Sigma(2^k,\mu2^l;m_1,m_2)=0$.
\end{itemize}
\noindent From now on assume that $\ell\geq 2$.
\begin{itemize}
\item If $k-\ell \geq 3$, then $\Sigma(2^k,\mu2^l;m_1,m_2)=$
\begin{multline*}2^{\ell}\sum_{\begin{smallmatrix}
0\leq i\leq \ell-2\\
i\equiv \ell\, (\text{mod }2)
\end{smallmatrix}}\big( g_{1}(2^k,m_1,2^{k-i})g_{1}(2^{\ell},m_2,2^{i+2})
+(-\mu)g_{-1}(2^k,m_1,2^{k-i})g_{-1}(2^{\ell},m_2,2^{i+2})\big).\end{multline*}
\item If $k-\ell=2$, then $\Sigma(2^k,\mu2^l;m_1,m_2)=$
\begin{multline*}(-1)^{\ell}2^{\ell}\sum_{\begin{smallmatrix}
0\leq i\leq \ell-2\\
i\equiv \ell\, (\text{mod }2)
\end{smallmatrix}}\big( g_{1}(2^k,m_1,2^{k-i})g_{1}(2^{\ell},m_2,2^{i+2})
+(-\mu)g_{-1}(2^k,m_1,2^{k-i})g_{-1}(2^{\ell},m_2,2^{i+2})\big).\end{multline*}
\item If $k-\ell=1$, then $\Sigma(2^k,\mu2^l;m_1,m_2)=$
\begin{multline*}(-\mu)^{\ell}2^{\ell}\sum_{\begin{smallmatrix}
0\leq i\leq \ell-2\\
i\equiv \ell\, (\text{mod }2)
\end{smallmatrix}}\big( g_{1}(2^{k},m_1,2^{k-i})g_{-1}(2^{\ell},m_2,2^{i+2})
+(-\mu)g_{-1}(2^{k},m_1,2^{k-i})g_{1}(2^{\ell},m_2,2^{i+2})\big).\end{multline*}
\end{itemize}


\end{proposition}
\noindent\textbf{Proof:} Proposition \ref{+pm2,k>l>0} implies that the set $S(2^k,\mu2^\ell)$ will be empty unless $\ell\geq2$. In this case the sum must be equal to $0$. From now on suppose that $\ell\geq 2$. 

We will begin with the formula for the splitting. Let $B_1=2^{i}b_1$, $B_2=2^jb_2$, and $\epsilon=\big(\frac{-1}{b_1}\big)$. Using equation (\ref{plucksplitting}), if $\gamma\in S(2^k,\mu2^l)$, then $s(\gamma) =\left(\frac{-\mu}{b_1}\right)
\left(\frac{b_1}{2^{k-\ell}}\right)\left(\frac{2^i}{2^{k-\ell}C_2+b_1b_2}\right)\left(\frac{2^{\ell-i}}{C_2}\right)$.
This computation will be broken into three cases, $k-\ell\geq3,$ $k-\ell=2$, $k-\ell=1$. 

First consider the case $k-\ell\geq 3$. By equation (\ref{plucksplitting}), $s(\gamma)=\left(\frac{-\mu}{b_1}\right)\left(\frac{b_1}{2^{k-\ell}}\right)\left(\frac{2^i}{b_1b_2}\right)\left(\frac{2^{\ell-i}}{C_2}\right)$.
Now we compute the exponential sum. By Proposition \ref{+pm2,k>l>0}, $\Sigma(2^k,\mu2^{\ell};m_1,m_2)=$
\begin{equation} \sum_{i=0}^{\ell-2}\sum_{b_1\in(\mathbb{Z}/2^{k-i}\mathbb{Z})^\times}\sum_{\begin{smallmatrix}
b_2\in(\mathbb{Z}/2^{i+2}\mathbb{Z})^\times\\
b_1b_2\equiv \mu \,(\text{mod } 4) 
\end{smallmatrix}}
\sum_{\begin{smallmatrix}
C_2\in (\mathbb{Z}/2^{\ell+2}\mathbb{Z})\\
C_2\equiv -1\,(\text{mod }4)
\end{smallmatrix}}\left(\frac{-\mu}{b_1}\right)\left(\frac{b_1}{2^{k+i-\ell}}\right)\left(\frac{2^i}{b_2}\right)\left(\frac{2^{\ell-i}}{C_2}\right)e^{2\pi i(\frac{m_1b1}{2^{k-i}}+\frac{\mu m_2b_2}{2^{i+2}})}.\label{2kleq1}
\end{equation}

First sum over $C_2$. If this is to be nonzero it must be that $\ell\equiv i \,(\text{mod } 2)$. Next consider the sum over $b_2$. By definition (\ref{GSumEven}) this is equal to $g_{\mu b_1}(2^{i},\mu m_2,2^{i+2})=g_{b_1}(2^{i},m_2,2^{i+2})$. Finally consider the sum over $b_1$. Begin by splitting this into two pieces based on the residue of $b_1$ mod $4$. By definition (\ref{GSumEven}), the summand in line (\ref{2kleq1}) corresponding to a fixed $i$ is equal to 
\begin{equation}2^{\ell}(g_{1}(2^k,m_1,2^{k-i})g_{1}(2^{i},m_2,2^{i+2})+(-\mu)g_{-1}(2^k,m_1,2^{k-i})g_{-1}(2^{i},m_2,2^{i+2})).\label{2kleq2}\end{equation}
Putting everything together yields
\begin{multline*}\Sigma(2^k,\mu2^l;m_1,m_2)=2^{\ell}\sum_{\begin{smallmatrix}
0\leq i\leq \ell-2\\
i\equiv \ell\, (\text{mod }2)
\end{smallmatrix}}\big( g_{1}(2^k,m_1,2^{k-i})g_{1}(2^{\ell},m_2,2^{i+2})\\
+(-\mu)g_{-1}(2^k,m_1,2^{k-i})g_{-1}(2^{\ell},m_2,2^{i+2})\big).\end{multline*}

Now we will consider the case $k-\ell=2$. By equation (\ref{plucksplitting}), $s(\gamma)=
\left(\frac{-\mu}{b_1}\right)\left(\frac{2^i}{4C_2+b_1b_2}\right)\left(\frac{2^{\ell-i}}{C_2}\right)$ $=(-1)^{i}\left(\frac{-\mu}{b_1}\right)\left(\frac{2^i}{b_1b_2}\right)\left(\frac{2^{\ell-i}}{C_2}\right)$.
By Proposition \ref{+pm2,k>l>0}, we have $\Sigma(2^k,-2^{\ell};m_1,m_2)=$
\begin{equation}\sum_{i=0}^{\ell-2}\sum_{b_1\in(\mathbb{Z}/2^{k-i}\mathbb{Z})^\times}\sum_{\begin{smallmatrix}
b_2\in(\mathbb{Z}/2^{i+2}\mathbb{Z})^\times\\
b_1b_2\equiv \mu \,(\text{mod } 4) 
\end{smallmatrix}}
\sum_{\begin{smallmatrix}
C_2\in (\mathbb{Z}/2^{\ell+2}\mathbb{Z})\\
C_2\equiv -1\,(\text{mod }4)
\end{smallmatrix}}(-1)^{i}\left(\frac{-\mu}{b_1}\right)\left(\frac{2^i}{b_1b_2}\right)\left(\frac{2^{\ell-i}}{C_2}\right)e^{2\pi i(\frac{m_1b1}{2^{k-i}}+\frac{\mu m_2b_2}{2^{i+2}})}.\label{aaa}
\end{equation}

First sum over $C_2$. If this is to be nonzero it must be that $\ell\equiv i \,(\text{mod } 2)$. Next consider the sum over $b_2$. By definition (\ref{GSumEven}) this is equal to $g_{\mu b_1}(2^{i},\mu m_2,2^{i+2})=g_{b_1}(2^{i},m_2,2^{i+2})$. Finally consider the sum over $b_1$. Begin by splitting this into two pieces based on the residue of $b_1$ mod $4$. By definition (\ref{GSumEven}) the sum in line (\ref{aaa}) associated with a fixed $i$ is equal to 
\[2^{\ell}(g_{1}(2^\ell,m_1,2^{k-i})g_{1}(2^{\ell},m_2,2^{i+2})+g_{-1}(2^\ell,m_1,2^{k-i})g_{-1}(2^{\ell},m_2,2^{i+2})).\]
Putting everything together and using that $k-\ell=2$ yields
\begin{multline*}\Sigma(2^k,-2^l;m_1,m_2)=(-1)^{\ell}2^{\ell}\sum_{\begin{smallmatrix}
0\leq i\leq \ell-2\\
i\equiv \ell\, (\text{mod }2)
\end{smallmatrix}}\big( g_{1}(2^k,m_1,2^{k-i})g_{1}(2^{\ell},m_2,2^{i+2})\\
+(-\mu)g_{-1}(2^k,m_1,2^{k-i})g_{-1}(2^{\ell},m_2,2^{i+2})\big).\end{multline*}

Lastly we will consider the case $k-\ell=1$. Note that since $b_1b_2\equiv -\mu \,(\text{mod }4)$ it follows that $\big(\frac{2^i}{6+b_1b_2}\big)=(-\mu)\big(\frac{2^i}{b_1b_2}\big)$. Thus, by equation (\ref{plucksplitting}), $s(\gamma)=
\left(\frac{-\mu}{b_1}\right)\left(\frac{b_1}{2}\right)\left(\frac{2^i}{2C_2+b_1b_2}\right)\left(\frac{2^{\ell-i}}{C_2}\right)=(-\mu)^i\left(\frac{-\mu}{b_1}\right)\left(\frac{b_1}{2}\right)\left(\frac{2^i}{b_1b_2}\right)\left(\frac{2^{\ell-i}}{C_2}\right)$.
Now we compute the exponential sum. By Proposition \ref{+pm2,k>l>0}, $\Sigma(2^k,-2^{\ell};m_1,m_2)=$
\begin{equation} \sum_{i=0}^{\ell-2}\sum_{b_1\in(\mathbb{Z}/2^{k-i}\mathbb{Z})^\times}\sum_{\begin{smallmatrix}
b_2\in(\mathbb{Z}/2^{i+2}\mathbb{Z})^\times\\
b_1b_2\equiv -\mu \,(\text{mod } 4) 
\end{smallmatrix}}\\
\sum_{\begin{smallmatrix}
C_2\in (\mathbb{Z}/2^{\ell+2}\mathbb{Z})\\
C_2\equiv -1\,(\text{mod }4)
\end{smallmatrix}}(-\mu)^{i}\left(\frac{-\mu}{b_1}\right)\left(\frac{b_1}{2}\right)\left(\frac{2^i}{b_1b_2}\right)\left(\frac{2^{\ell-i}}{C_2}\right)e^{2\pi i(\frac{m_1b1}{2^{k-i}}+\frac{\mu m_2b_2}{2^{i+2}})}.\label{aab}
\end{equation}

First sum over $C_2$. If this is to be nonzero it must be that $\ell\equiv i \,(\text{mod } 2)$. Next consider the sum over $b_2$. By definition (\ref{GSumEven}) this is equal to $g_{-\mu b_1}(2^{i},\mu m_2,2^{i+2})=g_{-b_1}(2^{i},m_2,2^{i+2})$. Finally consider the sum over $b_1$. Begin by splitting this into two pieces based on the residue of $b_1$ mod $4$. By definition (\ref{GSumEven}) the sum in line \ref{aab} corresponding to a fixed $i$ is equal to 
\[2^{\ell}((-\mu)^{\ell}g_{1}(2^{\ell+1},m_1,2^{k-i})g_{-1}(2^{\ell},m_2,2^{i+2})+(-\mu)^{\ell+1}g_{-1}(2^{\ell+1},m_1,2^{k-i})g_{1}(2^{\ell},m_2,2^{i+2})).\]
Putting everything together yields
\begin{multline*}\Sigma(2^k,-2^l;m_1,m_2)=(-\mu)^{\ell}2^{\ell}\sum_{\begin{smallmatrix}
0\leq i\leq \ell-2\\
i\equiv \ell\, (\text{mod }2)
\end{smallmatrix}}\big( g_{1}(2^{k},m_1,2^{k-i})g_{-1}(2^{\ell},m_2,2^{i+2})\\
+(-\mu)g_{-1}(2^{k},m_1,2^{k-i})g_{1}(2^{\ell},m_2,2^{i+2})\big).\end{multline*}
\EndProof


\begin{proposition}
\label{sigma +-2,k>0,m1,m2} Let $k\in \mathbb{Z}_{>0}$, and let $m_1,m_2\in\mathbb{Z}$. Then

\begin{multline*}
\Sigma(2^k,-2^{k};m_1,m_2)=2^{k}\delta_{2|k}(\delta_{2^k|m_1}2^k+\delta_{2^k|m_2}2^k-\delta_{2^{k-1}|m_1}2^{k-1}\\
+\sum_{i=1}^{k-1}(-1)^ig(2^2,m_1,2^{k-i})g(2^2,m_2,2^i)+\sum_{i=1 }^{k-1}\delta_{2^{i-1}|m_2}2^{i-1}g(2^2,m_1,2^{k-i})).
\end{multline*}
\end{proposition}

\textbf{Proof:} Recall the description of S$(2^k,-2^k)$ from Proposition \ref{+-2,k>0}. The computation of $\Sigma(2^k,-2^k;m_1,m_2)$ will utilize the cases discussed in that proposition.
We begin by computing the splitting. If $B_1$ and $B_2\neq0$, let $B_1=2^{i}b_1$, $B_2=2^jb_2$. By equation (\ref{plucksplitting}), 
$s(\gamma)=\left(\frac{2^i}{C_2+2^{2+i+j-k}b_1b_2}\right)\left(\frac{2^{k-i}}{C_2}\right)$. If $B_1$ or $B_2=0$, then $s(\gamma) = \left(\frac{2^k}{C_2}\right)$. 

In case 1 the sum is given by
\begin{equation*}
\sum_{\begin{smallmatrix}
C_2\in\mathbb{Z}/2^{2+k}\mathbb{Z}\\
C_2\equiv -1\,(\text{mod } 4)
\end{smallmatrix}}\left(\frac{2^k}{C_2}\right)=2^{k}\delta_{2|k}.
\end{equation*}
We will incorporate case 1 into both case 2 and case 3 and subtract the duplicated terms later. With the previous remark in mind, case 2 yields
\begin{equation*}
\sum_{B_2\in \mathbb{Z}/2^k\mathbb{Z}}\sum_{\begin{smallmatrix}
C_2\in\mathbb{Z}/2^{2+k}\mathbb{Z}\\
C_2\equiv -1\,(\text{mod } 4)
\end{smallmatrix}}\left(\frac{2^k}{C_2}\right)e^{2\pi i(-m_2\frac{B_2}{2^k})}=2^{2k}\delta_{2|k}\delta_{2^k|m_2}.
\end{equation*}

Similarly, for case 3 we have

\begin{equation*}
\sum_{B_1\in \mathbb{Z}/2^k\mathbb{Z}}\sum_{\begin{smallmatrix}
C_2\in\mathbb{Z}/2^{2+k}\mathbb{Z}\\
C_2\equiv -1\,(\text{mod } 4)
\end{smallmatrix}}\left(\frac{2^k}{C_2}\right)e^{2\pi i(m_1\frac{B_1}{2^k})}=2^{2k}\delta_{2|k}\delta_{2^k|m_1}.
\end{equation*}

Once we subtract the duplicate case 1 the sum of the first three cases is equal to

\begin{equation} \label{2kk-eq1}
2^{k}\delta_{2|k}(2^k\delta_{2^k|m_1}+2^k\delta_{2^k|m_2}-1).
\end{equation}

Finally we must consider case 4. The cases $i+j=k$ and $i+j>k$ must be dealt with separately. First consider $i+j=k$. In this case the sum is given by 
\begin{multline}
\sum_{i=1}^{k-1}\sum_{b_1\in (\mathbb{Z}/2^{k-i}\mathbb{Z})^\times}
\sum_{b_2\in (\mathbb{Z}/2^{i}\mathbb{Z})^\times}
\sum_{\begin{smallmatrix}
C_2\in\mathbb{Z}/2^{2+k}\mathbb{Z}\\
C_2\equiv -1\,(\text{mod } 4)
\end{smallmatrix}}(-1)^i\left(\frac{2^k}{C_2}\right)e^{2\pi i(m_1\frac{b_1}{2^{k-i}}- m_2\frac{b_2}{2^i})}\\
=2^k\delta_{2|k}\sum_{i=1}^{k-1}(-1)^ig(2^k,m_1,2^{k-i})g(2^k,- m_2,2^i).\label{2kk-eq2}
\end{multline}

Now suppose that $i+j>k$. In this case the sum is given by 
\begin{multline*}
\sum_{\begin{smallmatrix}
1\leq i\leq k-1\\
1\leq j\leq k-1\\
i+j>k
\end{smallmatrix}}
\sum_{b_1\in (\mathbb{Z}/2^{k-i}\mathbb{Z})^\times}
\sum_{b_2\in (\mathbb{Z}/2^{k-j}\mathbb{Z})^\times}
\sum_{\begin{smallmatrix}
C_2\in\mathbb{Z}/2^{2+k}\mathbb{Z}\\
C_2\equiv -1\,(\text{mod } 4)
\end{smallmatrix}}\left(\frac{2^k}{C_2}\right)e^{2\pi i(m_1\frac{b_1}{2^{k-i}}- m_2\frac{b_2}{2^{k-j}})}\\
=2^k\delta_{2|k}\sum_{\begin{smallmatrix}
1\leq i\leq k-1\\
1\leq j\leq k-1\\
i+j>k
\end{smallmatrix}}g(2^k,m_1,2^{k-i})g(2^k,-m_2,2^{k-j}).
\end{multline*}

This expression can be simplified using line (\ref{Gauss5}) in Lemma \ref{lemma:gauss}. Specifically,
\begin{equation*}
\sum_{1\leq k-j\leq i-1}g(2^2,m_2,2^{k-j}) = \delta_{2^{i-1}|m_2}2^{i-1}-1.
\end{equation*}
This can be summed over $i$ to get 
\begin{align}
\nonumber2^k\delta_{2|k}&\sum_{\begin{smallmatrix}
1\leq i\leq k-1\\
1\leq j\leq k-1\\
i+j>k
\end{smallmatrix}}g(2^k,m_1,2^{k-i})g(2^k,-m_2,2^{k-j})\\
\nonumber&=2^k\delta_{2|k}(\sum_{1\leq i \leq k-1}\delta_{2^{i-1}|m_2}2^{i-1}g(2^k,m_1,2^{k-i})-\sum_{1\leq k-i \leq k-1}g(2^k,m_1,2^{k-i}))\\
\label{2kk-eq3}&=2^k\delta_{2|k}(\sum_{1\leq i \leq k-1}\delta_{2^{i-1}|m_2}2^{i-1}g(2^k,m_1,2^{k-i})-(\delta_{2^{k-1}|m_1}2^{k-1}-1)).
\end{align}

Finally we add (\ref{2kk-eq1}), (\ref{2kk-eq2}), and (\ref{2kk-eq3}) together.\EndProof


\begin{proposition}
\label{sigma ++2,k>0,m1,m2} Let $k\in \mathbb{Z}_{>0}$, and let $m_1,m_2\in\mathbb{Z}$. Then
\begin{multline*}
\Sigma(2^k,2^k;m_1,m_2)=2^k\delta_{2|k}\sum_{i=0}^{k-1}(g_{1}(2^2,m_1,2^{k-i})+(-1)^{i+1}g_{-1}(2^2,m_1,2^{k-i}))\\
\times(g_{1}(2^{2},m_2,2^{i+1})+(-1)^ig_{-1}(2^{2},m_2,2^{i+1})).
\end{multline*}
\end{proposition}
\textbf{Proof:}
Recall the description of S$(2^k,2^k)$ from Proposition \ref{++2,k>0}. We begin by computing the splitting. Let $\gamma\in S(2^k,2^k)$. Let $B_1=2^{i}b_1$ and $B_2=2^jb_2$. By equation (\ref{plucksplitting}), we have $s(\gamma)=\big(\frac{-1}{b_1}\big)\left(\frac{2^i}{C_2+2b_1b_2}\right)\left(\frac{2^{k-i}}{C_2}\right)$.
Since $C_2\equiv -1 \,(\text{mod }4)$ and $(b_j,2)=1$, it follows that $C_2+2b_1b_2\equiv 1$ or $5$ $(\text{mod }8)$. Thus we will split $\Sigma(2^k,2^k;m_1,m_2)$ into two pieces based on the residue of $b_1b_2$ mod $4$. First assume that $b_1b_2\equiv 1\,\text{(mod } 4)$. In this case the sum is given by 

\begin{multline*}
\sum_{i=0}^{k-1}\sum_{b_1\in(\mathbb{Z}/2^{k-i}\mathbb{Z})^\times}\sum_{\begin{smallmatrix}
b_2\in(\mathbb{Z}/2^{i+1}\mathbb{Z})^\times\\
b_1b_2\equiv 1 \,(\text{mod } 4) 
\end{smallmatrix}}
\sum_{\begin{smallmatrix}
C_2\in (\mathbb{Z}/2^{\ell+2}\mathbb{Z})\\
C_2\equiv -1\,(\text{mod }4)
\end{smallmatrix}}\big(\frac{-1}{b_1}\big)\left(\frac{2^i}{C_2+2}\right)\left(\frac{2^{k-i}}{C_2}\right)e^{2\pi i(\frac{m_1b_1}{2^{k-i}}+\frac{m_2b_2}{2^{i+1}})}\\
=\sum_{i=0}^{k-1}\sum_{b_1\in(\mathbb{Z}/2^{k-i}\mathbb{Z})^\times}\sum_{\begin{smallmatrix}
b_2\in(\mathbb{Z}/2^{i+1}\mathbb{Z})^\times\\
b_1b_2\equiv 1 \,(\text{mod } 4) 
\end{smallmatrix}}
\sum_{\begin{smallmatrix}
C_2\in (\mathbb{Z}/2^{\ell+2}\mathbb{Z})\\
C_2\equiv -1\,(\text{mod }4)
\end{smallmatrix}}\big(\frac{-1}{b_1}\big)\left(\frac{2^{k}}{C_2}\right)e^{2\pi i(\frac{m_1b_1}{2^{k-i}}+\frac{m_2b_2}{2^{i+1}})}.
\end{multline*}

The sum over $C_2$ is nonzero if and only if $k$ is even. The sum over $b_2$ is equal to $g_{b_1}(2^2,m_2,2^{i+i})$. The sum over $b_1$ will be split into two cases based on the residue of $b_1$ mod $4$. Putting this together gives
\begin{multline*}
\sum_{i=0}^{k-1}\sum_{b_1\in(\mathbb{Z}/2^{k-i}\mathbb{Z})^\times}\sum_{\begin{smallmatrix}
b_2\in(\mathbb{Z}/2^{i+1}\mathbb{Z})^\times\\
b_1b_2\equiv 1 \,(\text{mod } 4) 
\end{smallmatrix}}
\sum_{\begin{smallmatrix}
C_2\in (\mathbb{Z}/2^{\ell+2}\mathbb{Z})\\
C_2\equiv -1\,(\text{mod }4)
\end{smallmatrix}}\big(\frac{-1}{b_1}\big)\left(\frac{2^i}{C_2+2}\right)\left(\frac{2^{k-i}}{C_2}\right)e^{2\pi i(\frac{m_1b_1}{2^{k-i}}+\frac{m_2b_2}{2^{i+1}})}\\
=\sum_{i=1}^{k-1}2^k\delta_{2|k}(g_{1}(2^{2},m_1,2^{k-i})g_{1}(2^{2},m_2,2^{i+1})-g_{-1}(2^{2},m_1,2^{k-i})g_{-1}(2^{2},m_2,2^{i+1})).
\end{multline*}

Next assume that $b_1b_2\equiv -1\,\text{(mod } 4)$. By similar computations we have  

\begin{multline*}
\sum_{i=0}^{k-1}\sum_{b_1\in(\mathbb{Z}/2^{k-i}\mathbb{Z})^\times}\sum_{\begin{smallmatrix}
b_2\in(\mathbb{Z}/2^{i+1}\mathbb{Z})^\times\\
b_1b_2\equiv -1 \,(\text{mod } 4) 
\end{smallmatrix}}
\sum_{\begin{smallmatrix}
C_2\in (\mathbb{Z}/2^{\ell+2}\mathbb{Z})\\
C_2\equiv -1\,(\text{mod }4)
\end{smallmatrix}}\big(\frac{-1}{b_1}\big)\left(\frac{2^i}{C_2-2}\right)\left(\frac{2^{k-i}}{C_2}\right)e^{2\pi i(\frac{m_1b_1}{2^{k-i}}+\frac{m_2b_2}{2^{i+1}})}\\
=\sum_{i=0}^{k-1}(-1)^i2^k\delta_{2|k}(g_{1}(2^{2},m_1,2^{k-i})g_{-1}(2^{2},m_2,2^{i+1})-g_{-1}(2^{2},m_1,2^{k-i})g_{1}(2^{2},m_2,2^{i+1})).
\end{multline*}

Putting everything together yields

\begin{multline*}
\Sigma(2^k,2^k;m_1,m_2)\\
= 2^k\delta_{2|k}\sum_{i=0}^{k-1}(g_{1}(2^{2},m_1,2^{k-i})g_{1}(2^{2},m_2,2^{i+1})-g_{-1}(2^{2},m_1,2^{k-i})g_{-1}(2^{2},m_2,2^{i+1})\\
+(-1)^ig_{1}(2^{2},m_1,2^{k-i})g_{-1}(2^{2},m_2,2^{i+1})-(-1)^ig_{-1}(2^{2},m_1,2^{k-i})g_{1}(2^{2},m_2,2^{i+1})).
\end{multline*}
Finally note that the summands can be factored.

\EndProof

\subsection{Big Cell Constant Term} \label{ssec:bigcellConstant}

We can specialize the results of the preceding subsections by setting $m_1=m_2=0$ to determine the big cell's contribution to the constant term Fourier coefficient. We assemble the complete constant term in Subsection \ref{ssec:Constant}.

\begin{proposition} \label{sigma+-p k,0,0,0} \label{sigma++,p,k>0,0,0} Let $p$ be an odd prime and let $k\in \mathbb{Z}$ such that $k>0$. 
Then:
\begin{itemize}
\item $\Sigma(1,-1;0,0)=1.$
\item$\Sigma(1,1;0,0)=0.$
\item$\Sigma(p^k,1;0,0)=0.$
\item If $k$ is odd, then $\Sigma(p^k,-1;0,0) = 0.$
\item If $k$ is even, then $\Sigma(p^k,-1;0,0) = \phi(p^k).$
\end{itemize}
\end{proposition}

\noindent\textbf{Proof:} This result follows directly from Proposition \ref{sigma+-p k,0,m1,m2} by setting $m_1=m_2=0$.\EndProof


\begin{proposition} \label{sigma+-p k>l>0,0,0}Let $p$ be an odd prime and let $k,l\in \mathbb{Z}$ such that $k>l>0.$ Then:
\begin{itemize}
\item $\Sigma(p^k,p^\ell;0,0)=0.$
\item If $p^{k-\ell}\equiv -1$ (mod $4$), then $\Sigma(p^k,-p^\ell;0,0)=0.$
\item If $p^{k-\ell}\equiv 1$ (mod $4$) and $k$ or $\ell$ is odd, then $\Sigma(p^k,-p^\ell;0,0)=0.$
\item If $p^{k-\ell}\equiv 1$ (mod $4$) and $k$ and $\ell$ are even, then 
\[\Sigma(p^k,-p^\ell;0,0)=\phi(p^{k-1})\phi(p^\ell)[(\frac{\ell+2}{2})p-(\frac{\ell-2}{2})].\]
\end{itemize}
\end{proposition}
\noindent\textbf{Proof:} This result follows directly from Proposition \ref{sigma+-p k>l>0,m1,m2} by setting $m_1=m_2=0$.\EndProof


\begin{proposition} \label{sigma+-p,k>0,0,0}\label{sigma++p,k>l>0,0,0}Let $p$ be an odd prime and let $k,l\in \mathbb{Z}$ such that $k>0.$ Then:
\begin{itemize}
\item $\Sigma(p^k,p^k;0,0)=0.$
\item If $k$ is odd, then $\Sigma(p^k,-p^k;0,0)=0.$
\item If $k$ is even, then
\end{itemize}
\[\Sigma(p^k,-p^k;0,0)=\phi(p^k)p^{k-2}((\frac{k+2}{2})p^2-(k-1)p+(\frac{k-2}{2})).\]

\end{proposition}
\noindent\textbf{Proof:} This result follows directly from Proposition \ref{sigma+-p,k>0,m1,m2} by setting $m_1=m_2=0$.\EndProof

\begin{proposition} \label{sigma 2,k>l>0,0,0} Let $k,\ell\in \mathbb{Z}$ such that $k>\ell\geq0$, and let $\mu=\pm1$. Then:
\begin{itemize}
\item If $\ell=0,1$, then $\Sigma(2^k,\mu2^l;0,0)=0.$
\item If $\ell\geq 2$, then $\Sigma(2^k,\mu2^l;0,0)=\delta_{k\equiv\ell\equiv0\,(\text{mod }2)}(\frac{1+(-\mu)}{2})2^{k+\ell-1}(\frac{\ell}{2}).$
\end{itemize}
\end{proposition}
\noindent\textbf{Proof:} This result follows directly from Proposition \ref{sigma 2,k>l>0,m1,m2} by setting $m_1=m_2=0$.\EndProof


\begin{proposition}
\label{sigma 2,k>0,0,0} Let $k\in \mathbb{Z}_{>0}$. Then $\Sigma(2^k,\mu2^{k};0,0)=\delta_{2|k}(\frac{1+(-\mu)}{2})(k+4)2^{2k-2}.$
\end{proposition}

\noindent\textbf{Proof:} This result follows directly from Proposition \ref{sigma +-2,k>0,m1,m2} and Proposition \ref{sigma ++2,k>0,m1,m2} by setting $m_1=m_2=0$.\EndProof


\subsection{The Formula of Brubaker-Bump-Friedberg-Hoffstein} \label{sec:BBFHform}


In this section the formulas for $\Sigma(p^k,\pm p^\ell,m_1,m_2)$, computed in Section \ref{ssec:Sigma2}, are manipulated to resemble those contained in Brubaker-Bump-Friedberg-Hoffstein \cite{BBFH07}.

\begin{proposition} Let $p$ be an odd prime, let $\mu=\pm1$, and let $k,\ell,r_1,r_2\in \mathbb{Z}$ such that $k>\ell>0$, $p^{k-\ell}\equiv -\mu\,(\text{mod }4)$, and $r_1,r_2\geq0.$ Then
\begin{align}\Sigma(p^k,\mu p^\ell;p^{r_1},p^{r_2})
=&\sum_{i=0}^{\ell}g(p^i,p^i,p^i)g(p^{\ell-i}, p^{r_2},p^{\ell-i})g(p^k,p^{r_1+\ell-i},p^k). \label{expeqn5}
\end{align}
\end{proposition}

\noindent\textbf{Proof:}  In Proposition \ref{sigma+-p k>l>0,m1,m2} it was shown that 
\begin{align}
\nonumber\Sigma(p^k,\mu p^\ell;p^{r_1}, p^{r_2})=&p^\ell g(p^\ell, p^{r_2},p^\ell)g(p^k,p^{r_1},p^{k-\ell})\\
&+\phi(p^\ell)\sum_{0\leq i< \ell}(\frac{1+(-1)^{\ell-i}}{2})g(p^i, p^{r_2},p^i)g(p^k,p^{r_1},p^{k-i}).\label{expeqn4}
\end{align}
First note that in line (\ref{expeqn5}) the terms with $l\equiv i+1$ (mod 2) vanish. Now the $i$th term on the right hand side of line (\ref{expeqn5}) will be paired with the $(l-i)$th term in the right hand side of line (\ref{expeqn4}), where $p^lg(p^l,p^{r_2},p^l)g(p^k,p^{r_1},p^{k-l})$ is considered to be the $i=l$ term. Now compute the terms explicitly using Lemma \ref{lemma:gauss} and note that they agree.\EndProof

\begin{proposition}  Let $p$ be an odd prime and let $k,r_1,r_2\in \mathbb{Z}$ such that $k>0$ and $r_1,r_2\geq0.$ Then
\begin{multline*}
\Sigma(p^k,-p^k;p^{r_1}, p^{r_2})
=\sum_{0\leq i\leq k-1}g(p^i,- p^{r_2},p^i)g(p^k,p^{r_1+i},p^k)g(p^{k-i}, p^{r_2+k-2i},p^{k-i})\\
+\delta_{k\leq r_2}p^kg(p^k,p^k,p^k).
\end{multline*}
\end{proposition}

\noindent\textbf{Proof:}  First suppose that $k$ is odd. In this case the claim is that 
\[\sum_{i=1}^{k-1}g(p^i,- p^{r_2},p^i)p^{k-i}g(p^{k-i}, p^{r_2},p^i)g(p^k,p^{r_1},p^{k-i})\]
\[= \sum_{0\leq i\leq k-1}g(p^i,- p^{r_2},p^i)g(p^k,p^{r_1+i},p^k)g(p^{k-i}, p^{r_2+k-2i},p^{k-i})+\delta_{k\leq r_2}p^kg(p^k,p^k,p^k),\]
by the computation of Proposition \ref{sigma+-p,k>0,m1,m2}.
By Lemma \ref{lemma:gauss},
\begin{align*}&g(p^i,- p^{r_2},p^i)p^{k-i}g(p^{k-i}, p^{r_2},p^i)g(p^k,p^{r_1},p^{k-i})\\
=&g(p^i,- p^{r_2},p^i)p^{k-2i}g(p^{k-i}, p^{r_2},p^i)g(p^k,p^{r_1+i},p^{k})\\
=&g(p^i,- p^{r_2},p^i)g(p^{k-i}, p^{r_2+k-2i},p^{k-i})g(p^k,p^{r_1+i},p^{k}).\end{align*}
Note that $k\equiv 1$ (mod $2$) was not needed for this. Thus it remains to prove that 
\[0= g(p^0,- p^{r_2},p^0)g(p^k,p^{r_1},p^k)g(p^k, p^{r_2+k},p^k)+\delta_{k\leq r_2}p^kg(p^k,p^k,p^k).\]
As $k$ is odd $\delta_{k\leq r_2}g(p^k,p^k,p^k)=0$ and $g(p^k, p^{r_2+k},p^k)=0$. Thus the desired equality holds.

Now suppose that $k$ is even. In this case it must be shown that 
\begin{align*}\sum_{i=1}^{k-1}&g(p^i,- p^{r_2},p^i)p^{k-i}g(p^{k-i}, p^{r_2},p^i)g(p^k,p^{r_1},p^{k-i})+\phi(p^k)(\delta_{k\leq r_1}p^k+\delta_{k\leq r_2}p^k-\delta_{k-1\leq r_1}p^{k-1})\\
=& \sum_{0\leq i\leq k-1}g(p^i,- p^{r_2},p^i)g(p^k,p^{r_1+i},p^k)g(p^{k-i}, p^{r_2+k-2i},p^{k-i})+\delta_{k\leq r_2}p^kg(p^k,p^k,p^k),
\end{align*} again by Proposition \ref{sigma+-p,k>0,m1,m2}.
As was mentioned in the previous case 
\begin{multline*}g(p^i,- p^{r_2},p^i)p^{k-i}g(p^{k-i}, p^{r_2},p^i)g(p^k,p^{r_1},p^{k-i})\\
=g(p^i,- p^{r_2},p^i)g(p^{k-i}, p^{r_2+k-2i},p^{k-i})g(p^k,p^{r_1+i},p^{k}).\end{multline*}
First note that since $k$ is even $\delta_{k\leq r_2}\phi(p^k)p^k=\delta_{k\leq r_2}p^kg(p^k,p^k,p^k)$. Thus it remains to show that $\phi(p^k)(\delta_{k\leq r_1}p^k-\delta_{k-1\leq r_1}p^{k-1})
=g(p^0,- p^{r_2},p^0)g(p^k,p^{r_1},p^k)g(p^k, p^{r_2+k},p^k)$.

Now $g(p^0,- p^{r_2},p^0) = 1$ and as $k$ is even $g(p^k, p^{r_2+k},p^k) = \phi(p^k)$. Thus it remains to show that $(\delta_{k\leq r_1}p^k-\delta_{k-1\leq r_1}p^{k-1})
=g(p^k,p^{r_1},p^k)$. This can be proved directly by considering the cases $r_1<k-1,$ $r_1=k-1$, and $r_1\geq k$.
\EndProof




\section{Metaplectic Eisenstein Distribution}\label{sec:MetaEisenstein}


Recall the definition of the metaplectic Eisenstein distribution 
\begin{equation}\label{metaEisen2}\tilde{E}(\tilde{g},\lambda)
 = \sum_{\gamma\in \Gamma_{\infty}\backslash\Gamma}\pi(S(\gamma)^{-1})\tilde{\tau}( \tilde{g})\in \widetilde{V}^{-\infty}_{\lambda,\phi}\end{equation}
first intoduced on line (\ref{metaEisen}). In this section we compute the Whittaker distributions of $\tilde{E}(\tilde{g},\lambda)$. 

The Eisenstein series may be split into six parts depending on the Bruhat cell in which $\gamma$ resides.  Let
\begin{equation}\label{metaEisenw}\tilde{E}_{w}(\tilde{g},\lambda) = \sum_{\gamma\in \Gamma_{\infty}\backslash\Gamma\cap (NwB)}
\pi(\tilde{\gamma}^{-1})\tilde{\tau}(\tilde{g}).\end{equation}
Thus, $\tilde{E}(\tilde{g},\lambda) = \sum_{w\in W}\tilde{E}_w(\tilde{g},\lambda).$ As usual, the Fourier coefficient computation is accomplished by addressing each $\tilde{E}_w$ individually.

The subsections \ref{ssec:tildeB}-\ref{ssec:tildeNwlB} contain the formulas for the Fourier coefficients. For each Bruhat cell the computations are similar, so we will only include the details for $B$, $Nw_{\alpha_{1}}w_{\alpha_{2}}B$, and $Nw_{\ell}B$. Let $\psi_{m_1,m_2}(
\left(\begin{smallmatrix}
1 & x & z\\
0 & 1 & y\\
0 & 0 & 1 \end{smallmatrix}\right))
=e^{2\pi i(m_1x+m_2y)}$. In the computations that follow, $f\in \widetilde{V}_{-\lambda,^\intercal\phi^{-1}}^{\infty}$ and thus may be paired against elements of $\widetilde{V}_{\lambda,\phi}^{-\infty}$ as described in Section \ref{ssec:PSeries}. In particular, $f$ may be paired against $\tilde{E}$.


The Whittaker distribution of $\tilde{E}$ is defined to be
\begin{equation}\label{metaWhit}\int_{\Gamma_{\infty}\backslash N}\tilde{E}(n\tilde{g},\lambda)\psi_{m_1,m_2}(n)dn.\end{equation}
In what follows the parameters $m_1$, $m_2$, and $\lambda$ will be suppressed.


\subsection{Bruhat Cell: $B$} \label{ssec:tildeB}

\begin{proposition} \label{metawhitid}If $f=
\left[\begin{smallmatrix}
f_1\\
f_2
\end{smallmatrix}\right]\in\widetilde{V}_{-\lambda,^\intercal\phi^{-1}}^{\infty}$, then
\begin{equation} \langle f(\tilde{g}),\int_{\Gamma_{\infty}\backslash N}\tilde{E}_{\text{id}}(n\tilde{g})\psi(n)dn\rangle_{\lambda,\phi}=
\int_{\Gamma_{\infty}\backslash N}\psi(n)dn
f_1((w_{\ell},1)).
\end{equation}
This distribution is nonzero if and only if $m_1=m_2=0$.
\end{proposition}

\noindent\text{Proof:} This follows from the definition of $\tilde{\tau}$ and the left $N$-invariance of $\tilde{E}_{\text{id}}(\tilde{g})=\tilde{\tau}(\tilde{g})$.\EndProof

\subsection{Bruhat Cell: $Nw_{\alpha_{1}}B$} \label{tildeNw1B}



Let $w=\left( \left(\begin{smallmatrix}
 & 1 & \\
 &  & -1\\
-1 &  &  \end{smallmatrix}\right),1 \right)=\left( \left(\begin{smallmatrix}
  & -1 & \\
-1 &  & \\
&  & -1 \end{smallmatrix} \right)w_{\ell},1\right).$ We begin by computing the effect of each summand of $\widetilde{E}_{w_{\alpha_1}}$ on a test vector.

\begin{proposition}\label{distalpha1} Let $\gamma\in Nw_{\alpha_1}B$ with Pl\"{u}cker coordinates $(0,0,-1,0,4B_2,C_2)$ and let $f=
\left[\begin{smallmatrix}
f_1\\
f_2
\end{smallmatrix}\right]\in \widetilde{V}_{-\lambda,^\intercal\phi^{-1}}^{\infty}$. Then:\\
\[ \langle f,\pi(\gamma^{-1})\tilde{\tau}\rangle_{\lambda,\phi}= \left(\frac{B_2}{-C_2}\right)|4B_2|^{-1-\lambda_2+\lambda_3}
f_2(wn(0,\frac{-C_2 }{4B_2 },0))
\begin{cases}
-i,&\text{ if } B_2>0;\\
1,& \text{ if } B_2<0.
\end{cases}\]

\end{proposition}


Now we compute the Fourier coefficients of $\widetilde{E}_{w_{\alpha_1}}$.

\begin{proposition}\label{metawhitalpha1}
Let $f=
\left[\begin{smallmatrix}
f_1\\
f_2
\end{smallmatrix}\right]\in \widetilde{V}_{-\lambda,^\intercal\phi^{-1}}^{\infty}$.

If $m_2\neq0$, then $\langle f(\tilde{g}),\int_{\Gamma_{\infty}\backslash N}\tilde{E}_{w_{\alpha_1}}(n\tilde{g})\psi(n)dn\rangle_{\lambda,\phi}=0.$

If $m_2=0$, then
\begin{multline} 
\langle f(\tilde{g}),\int_{\Gamma_{\infty}\backslash N}\tilde{E}_{w_{\alpha_1}}(n\tilde{g})\psi(n)dn\rangle_{\lambda,\phi}
=\left(\sum_{B_{2}\in\mathbb{Z}_{>0}}(4B_2)^{-1-\lambda_2+\lambda_3}K_{1}(m_1;4B_2)\right)\\
\times\int_{\mathbb{R}}f_2(n(x,0,0)w)e^{-2\pi im_1x}dx.
\end{multline}
\end{proposition}


This L-function is studied by Shimura \cite{S75} and Bate \cite{B}. Bate introduces the L-function in Section 4 and provides an explicit description of it in terms of quadratic L-functions in Proposition 5.2.

\subsection{Bruhat Cell: $Nw_{\alpha_{2}}B$} \label{ssec:tildeNw2B}


Let $w=\left(\left( \begin{smallmatrix}
 &  & -1\\
-1&  & \\
& 1 &  \end{smallmatrix} \right), 1\right)
=\left(\left( \begin{smallmatrix}
-1 &  & \\
&  & -1\\
& -1 &  \end{smallmatrix} \right)
w_{\ell}, 1\right)$. We begin by computing the effect of each summand of $\widetilde{E}_{w_{\alpha_2}}$ on a test vector.

\begin{proposition}\label{distalpha2} Let $\gamma\in Nw_{\alpha_2}B$ with Pl\"{u}cker coordinates $(0,4B_1,C_1,0,0,-1)$ and let $f=
\left[\begin{smallmatrix}
f_1\\
f_2
\end{smallmatrix}\right]\in \widetilde{V}_{-\lambda,^\intercal\phi^{-1}}^{\infty}$. Then:\\
\[\langle f,\pi(\gamma^{-1})\tilde{\tau}\rangle_{\lambda,\phi}=\left(\frac{-B_1}{-C_1}\right)|4B_1|^{-1-\lambda_1+\lambda_2}
\begin{cases}\phantom{-}if_1(wn(\frac{C_1}{4B_1},0,0)),&\text{ if } B_{1}>0;\\
-\phantom{i}f_2(wn(\frac{C_1}{4B_1},0,0)),&\text{ if } B_{1}<0.
\end{cases}\]
\end{proposition}


Now we compute the Fourier coefficients of $\widetilde{E}_{w_{\alpha_2}}$. 
 
 \begin{proposition}\label{metawhitalpha2}
 Let $f=
\left[\begin{smallmatrix}
f_1\\
f_2
\end{smallmatrix}\right]\in \widetilde{V}_{-\lambda,^\intercal\phi^{-1}}^{\infty}$.

If $m_1\neq0$, then $\langle f(\tilde{g}),\int_{\Gamma_{\infty}\backslash N}\tilde{E}_{w_{\alpha_2}}(n\tilde{g})\psi(n)dn\rangle_{\lambda,\phi}=0.$

If $m_1=0$, then
\begin{align*} 
&\langle f(\tilde{g}),\int_{\Gamma_{\infty}\backslash N}\tilde{E}_{w_{\alpha_2}}(n\tilde{g})\psi(n)dn\rangle_{\lambda,\phi}\\
=&i\left(\sum_{B_1\in\mathbb{Z}_{>0}}
(4B_1)^{-1-\lambda_1+\lambda_2}\left(\frac{K_{1}(m_2;4B_1)+K_{-1}(m_2;4B_1)}{2}\right)\right)\\
&\hspace{7cm}\times\int_{\mathbb{R}}f_1\left(n(0,y,0)
\left(w,1\right)
\right)\psi^{-1}(n(0,y,0))dy\\
&+(-i)\left(\sum_{B_1\in\mathbb{Z}_{>0}}
(4B_1)^{-1-\lambda_1+\lambda_2}
\left(\frac{K_{1}(m_2;4B_1)-K_{-1}(m_2;4B_1)}{2}\right)\right)\\
&\hspace{7cm}\times\int_{\mathbb{R}}f_2\left(n(0,y,0)
\left(w,1\right)
\right)\psi^{-1}(n(0,y,0))dy.\\
\end{align*}
\end{proposition}


\subsection{Bruhat Cell: $Nw_{\alpha_{1}}w_{\alpha_{2}}B$} \label{ssec:tildeNw1w2B}


Let $w=\left(\left( \begin{smallmatrix}
 & -1 & \\
1&  & \\
&  & 1 \end{smallmatrix} \right), 1\right)
=\left(\left( \begin{smallmatrix}
 & 1 & \\
&  & 1\\
1&  &  \end{smallmatrix} \right)
w_{\ell}, 1\right)$. We begin by computing the effect of each summand of $\widetilde{E}_{w_{\alpha_1}w_{\alpha_2}}$ on a test vector.

\begin{proposition}\label{distmetaalpha1alpha2} Let $\gamma\in Nw_{\alpha_1}w_{\alpha_2}B$ with Pl\"{u}cker coordinates $(0,4B_1,C_1,4A_2,4B_2,C_2)$ and $f=
\left[\begin{smallmatrix}
f_1\\
f_2
\end{smallmatrix}\right]\in \widetilde{V}_{-\lambda,^\intercal\phi^{-1}}^{\infty}$, then:
\begin{multline*}\langle f,\pi(\gamma^{-1})\tilde{\tau}\rangle_{\lambda,\phi}=|4B_1|^{-1-\lambda_1+\lambda_2}|4A_2|^{-1-\lambda_2+\lambda_3}\left(\frac{A_2/B_1}{-C_2}\right)\left(\frac{-B_1}{-C_1}\right)\\
\times\begin{cases}
\phantom{i}f_2\left(wn(0,\frac{C_2}{4A_2},\frac{B_2}{A_2})\right),&\text{ if } B_1>0,A_2>0;\\
if_2\left(wn(0,\frac{C_2}{4A_2},\frac{B_2}{A_2})\right),&\text{ if } B_1>0,A_2<0;\\
\phantom{i}f_1\left(wn(0,\frac{C_2}{4A_2},\frac{B_2}{A_2})\right),&\text{ if } B_1<0,A_2<0;\\
if_1\left(wn(0,\frac{C_2}{4A_2},\frac{B_2}{A_2})\right),&\text{ if } B_1<0,A_2>0.
\end{cases}
\end{multline*}

\end{proposition}

\textbf{Proof:} By the definition of $\tilde{\tau}$ we have
\begin{equation}
\langle f,\pi(\gamma^{-1})\tilde{\tau}\rangle_{\lambda,\phi}=\langle \pi(\gamma)f,\tilde{\tau}\rangle_{\lambda,\phi}=\int_{\mathbb{R}^3}(\pi(\gamma)f)\left((w_{\ell},1)\left( \begin{smallmatrix}
1 & x & z \\
0 & 1 & y \\
0 & 0 & 1 \end{smallmatrix} \right)\right)\tilde{\tau}\left((w_{\ell},1)\left( \begin{smallmatrix}
1 & x & z \\
0 & 1 & y \\
0 & 0 & 1 \end{smallmatrix} \right)\right)dxdydz.\label{TempEqn1}
\end{equation}
By the definition of $\tilde{\tau}$, line (\ref{TempEqn1}) is equal to the first component of the $f((\gamma,s(\gamma))^{-1}(w_{\ell},1))$. We can compute this quantity using equation (\ref{walpha12iden}), the formula for $s(\gamma)$ from Proposition \ref{PluckSplitAll}, and the definition of $\phi$ from equation (\ref{Q8rep}).\EndProof

Now we compute the Fourier coefficients of $\widetilde{E}_{w_{\alpha_1}w_{\alpha_2}}$.

\begin{proposition}\label{metawhitalpha12}
Let $f=
\left[\begin{smallmatrix}
f_1\\
f_2
\end{smallmatrix}\right]\in \widetilde{V}_{-\lambda,^\intercal\phi^{-1}}^{\infty}$.

If $m_1\neq0$, then $\langle f(\tilde{g}),\int_{\Gamma_{\infty}\backslash N}\tilde{E}_{w_{\alpha_1}w_{\alpha_2}}(n\tilde{g})\psi(n)dn\rangle_{\lambda,\phi}=0.$

If $m_1=0$, then
\begin{align*} 
&\langle f(\tilde{g}),\int_{\Gamma_{\infty}\backslash N}\tilde{E}_{w_{\alpha_1}w_{\alpha_2}}(n\tilde{g})\psi(n)dn\rangle_{\lambda,\phi}\\
=&(-1+i)\left(\sum_{
B_1\in\mathbb{Z}_{>0}}|4B_1|^{-1-\lambda_1+\lambda_3}
\left(\frac{K_{1}(m_2;4B_1)-K_{-1}(m_2;4B_1)}{2}\right)\right)\\
&\hspace{3cm}\times\left(2^{-2(-1-\lambda_2+\lambda_3)}\frac{\zeta(2\lambda_2-2\lambda_3)}{\zeta_2(2\lambda_2-2\lambda_3+1)}\right)
\int_{\mathbb{R}^2}f_1\left(n(0,y,z)w\right)e^{-2\pi im_2y}dydz\\
&+(1+i)\left(\sum_{
B_1\in\mathbb{Z}_{>0}}|4B_1|^{-1-\lambda_1+\lambda_3}
\left(\frac{K_{1}(m_2;4B_1)+K_{-1}(m_2;4B_1)}{2}\right)\right)\\
&\hspace{3cm}\times\left(2^{-2(-1-\lambda_2+\lambda_3)}\frac{\zeta(2\lambda_2-2\lambda_3)}{\zeta_2(2\lambda_2-2\lambda_3+1)}\right)\int_{\mathbb{R}^2}f_2\left(n(0,y,z)w\right)e^{-2\pi im_2y}dydz.
\end{align*}
\end{proposition}


\textbf{Proof:} We begin with a bit of notation. In what follows, the summation over $\gamma$ (as opposed to over $\gamma^\prime$) will be indexed by P\"{u}cker coordinates for the set $\Gamma_{\infty}\backslash \Gamma$ such that the matrix representative $\gamma$ is in $Nw_{\alpha_1}w_{\alpha_2}B$. In terms of the Pl\"{u}cker coordinates, this set can be described as 
\begin{multline}\label{12set}S_{\alpha_1\alpha_2}=\{(0,4B_1,C_1,4A_2,4B_2,C_2)\in\mathbb{Z}^{6}| B_1,A_2\neq0,(B_1,C_1)=1,\\4B_1 B_2=-C_1 A_2,\,(\frac{A_2}{B_1},C_2)=1,\,C_{j}\equiv -1\text{ (mod }4)\}.\end{multline}
The summation over $\gamma^\prime$ will consist of distinct representatives of the double coset space $\Gamma_{\infty}\backslash \Gamma/(\Gamma_{\infty}\cap w\Gamma_{\infty}w^{-1})$. Proposition \ref{prop:PluckerSym} shows that this double coset space is in bijection with the subset of $S_{\alpha_1\alpha_2}$ such that $0\leq C_1<|4B_1|$ and $0\leq C_2  <|4A_2|$. The switch between $\gamma$ and $\gamma^\prime$ occurs when the integral over $\Gamma_{\infty}\backslash N$ is unfolded.

Let $n$ denote $n(x,y,z)$, and define $\epsilon_{\gamma}\in\{\pm1,\pm i\}$ and $f_{\gamma}=f_{j}$,with $j\in\{1,2\}$, such that 
\[\langle \pi(n)f,\pi(S(\gamma))^{-1}\tilde{\tau}\rangle_{\lambda,\phi}=|4B_1|^{-1-\lambda_1+\lambda_2}|4A_2|^{-1-\lambda_2+\lambda_3}\epsilon_{\gamma}s(\gamma)f_{\gamma}\left( n^{-1}wn(0,\frac{C_2^\prime }{A_2^\prime},\frac{B_2^\prime }{A_2^\prime})\right),\]
in accordance with Proposition \ref{distmetaalpha1alpha2}. Observe that $\epsilon_{\gamma}$, $s(\gamma)$, and $j$ only depend on the double coset $\Gamma_{\infty}\backslash \Gamma/(\Gamma_{\infty}\cap w\Gamma_{\infty}w^{-1})$.
By Proposition \ref{distmetaalpha1alpha2},
\begin{align}\nonumber\langle f(\tilde{g}),\int_{\Gamma_{\infty}\backslash N}&\tilde{E}_{w_{\alpha_1}w_{\alpha_2}}(n\tilde{g})\psi(n)dn\rangle_{\lambda,\phi}=\int_{\Gamma_{\infty}\backslash N}\sum_{
\gamma\in \Gamma_\infty \backslash\Gamma\cap(
Nw_{\alpha_1}w_{\alpha_2}B)}\langle f,\pi((\gamma n)^{-1})\tilde{\tau}\rangle_{\lambda,\phi}\,\psi(n)dn\\
=&\int_{\Gamma_{\infty}\backslash N}\sum_{\gamma}|4B_1|^{-1-\lambda_1+\lambda_2}|4A_2|^{-1-\lambda_2+\lambda_3} \epsilon_{\gamma}s(\gamma)f_{\gamma}\left(n^{-1}wn(0,\frac{C_2}{4A_2},\frac{B_2 }{A_2})\right)\psi\left(n\right)dn\label{SL3Rwhit12eqn1}.
\end{align}
Since the $2$-cocycle $\sigma$ is trivial when one entry is an element of $N$, we have $w n(0,\frac{C_2^\prime }{A_2^\prime},\frac{-C_1^\prime}{B_1^\prime})=n(0,\frac{B_2^\prime }{A_2^\prime},\frac{-C_2^\prime }{A_2^\prime})w$. Thus,  
\begin{equation}
(\ref{SL3Rwhit12eqn1})=\int_{\Gamma_{\infty}\backslash N}\sum_{\gamma}|4B_1|^{-1-\lambda_1+\lambda_2}|4A_2|^{-1-\lambda_2+\lambda_3} \epsilon_{\gamma}s(\gamma)f_{\gamma}\left((n(0,\frac{-B_2 }{A_2},\frac{C_2 }{4A_2})n)^{-1}w\right)\psi\left(
n\right)dn.\label{SL3Rwhit12eqn2}
\end{equation}
The next step requires unfolding the integral. An element $\gamma\in Nw_{\alpha_1}w_{\alpha_2}B\cap\Gamma$ that represents a coset in $\Gamma_{\infty}\backslash \text{SL}(3,\mathbb{Z})$ can be factored as $\gamma=\gamma^{\prime}\gamma^{\prime\prime}$, where $\gamma^\prime\in Nw_{\alpha_1}w_{\alpha_2}B\cap\Gamma$ represents a double coset in $\Gamma_{\infty}\backslash \Gamma/(\Gamma_{\infty}\cap w\Gamma_{\infty}w^{-1})$ and $\gamma^{\prime\prime}\in(\Gamma_{\infty}\cap w\Gamma_{\infty}w^{-1})\cong\mathbb{Z}^2$. If a set of distinct representatives of the double coset space is identified, then the factorization is unique. The integral over $\Gamma_{\infty}\backslash N$ is unfolded with respect to the sum over $\Gamma_{\infty}\cap w\Gamma_{\infty}w^{-1}$; the result is an integral over $(\Gamma_{\infty}\cap w\Gamma_{\infty}^{op}w^{-1})\backslash N$. Thus
\begin{align}\nonumber(\ref{SL3Rwhit12eqn2})=&\int_{\Gamma_{\infty}\backslash N}\sum_{\gamma^\prime}\sum_{j,k\in\mathbb{Z}}|4B_1|^{-1-\lambda_1+\lambda_2}|4A_2|^{-1-\lambda_2+\lambda_3}\epsilon_{\gamma^\prime}s(\gamma^\prime)\\
\nonumber&\hspace{1.5cm}\times f_{\gamma^\prime}\left((n(0,\frac{-B_2 }{A_2},\frac{C_2}{4A_2})n(x,y+j,z+k))^{-1}w\right)\psi\left(
n(x,y+j,z+k)\right)dn\\
\nonumber=&\sum_{\gamma^\prime}|4B_1|^{-1-\lambda_1+\lambda_2}|4A_2|^{-1-\lambda_2+\lambda_3}\epsilon_{\gamma^\prime}s(\gamma^\prime)\\
&\hspace{1.5cm}\times \int_{(\Gamma_{\infty}\cap w\Gamma_{\infty}^{op}w^{-1})\backslash N}f_{\gamma^\prime}\left((n(0,\frac{-B_2 }{A_2},\frac{C_2}{4A_2})n(x,y,z))^{-1}w\right)\psi\left(
n(x,y,z)\right)dn.\label{SL3Rwhit12eqn3}
\end{align}
Perform the change of variables $n\mapsto nn(0,\frac{B_2 }{A_2},\frac{-C_2}{4A_2})$ and recall that $\frac{B_2 }{A_2}=\frac{-C_1}{4B_1}$ to see that 
\begin{multline}(\ref{SL3Rwhit12eqn3})=\sum_{\gamma^\prime}|4B_1|^{-1-\lambda_1+\lambda_2}|4A_2|^{-1-\lambda_2+\lambda_3}\epsilon_{\gamma^\prime}s(\gamma^\prime)e^{2\pi i(m_2\frac{-C_1}{4B_1})}\\
\times\int_{(\Gamma_{\infty}\cap w\Gamma_{\infty}^{op}w^{-1})\backslash N}f_{\gamma^\prime}\left((n(0,\frac{-B_2 }{A_2},\frac{C_2}{4A_2})nn(0,\frac{B_2 }{A_2},\frac{-C_2}{4A_2}))^{-1}w\right)\psi\left(
n\right)dn.\label{SL3Rwhit12eqn4}
\end{multline}
The previous change of variables is advantageous as elements of the form $n n^\prime n^{-1}$ in $N$ can be written in the form $n^\prime z$, where $z\in[N,N]$ is an element in the derived subgroup. In this case, $n(0,\frac{-B_2^\prime }{A_2^\prime},\frac{C_2^\prime}{A_2^\prime})n(x,y,z)n(0,\frac{B_2^\prime }{A_2^\prime},\frac{-C_2^\prime}{A_2^\prime})=n(x,y,z)n(0,0,\frac{-B_2^\prime }{A_2^\prime}x)$. Thus another change of variables $n(x,y,z)\mapsto n(x,y,z)n(0,0,\frac{B_2^\prime }{A_2^\prime}x)$ pushes the element of the derived subgroup into the character where it contributes trivially. Thus
\begin{align}(\ref{SL3Rwhit12eqn4})=&\sum_{\gamma^\prime}|4B_1|^{-1-\lambda_1+\lambda_2}|4A_2|^{-1-\lambda_2+\lambda_3}\epsilon_{\gamma^\prime}s(\gamma^\prime)e^{2\pi i(m_2\frac{-C_1}{4B_1})}\\
\nonumber&\hspace{1cm}\times\int_{(\Gamma_{\infty}\cap w\Gamma_{\infty}^{op}w^{-1})\backslash N}
f_{\gamma^\prime}\left((n(x,y,z)n(0,0,\frac{-B_2^\prime }{A_2^\prime}x))^{-1}w\right)\psi\left(
n(x,y,z)\right)dn\\
=\sum_{\gamma^\prime}&|4B_1|^{-1-\lambda_1+\lambda_2}|4A_2|^{-1-\lambda_2+\lambda_3}\epsilon_{\gamma^\prime}s(\gamma^\prime)e^{2\pi i(m_2\frac{-C_1}{4B_1})}\int_{(\Gamma_{\infty}\cap w\Gamma_{\infty}^{op}w^{-1})\backslash N}f_{\gamma^\prime}\left(n^{-1}w\right)\psi\left(n\right)dn.\label{SL3Rwhit12eqn5}
\end{align}
We now perform the change of variables $n\mapsto n^{-1}$ and apply the right $\widetilde{N}_{-}$-invariance to remove the $x$-variable to see that
\begin{multline}(\ref{SL3Rwhit12eqn5})=\sum_{\gamma^\prime}|4B_1|^{-1-\lambda_1+\lambda_2}|4A_2|^{-1-\lambda_2+\lambda_3}\epsilon_{\gamma^\prime}s(\gamma^\prime)e^{2\pi i(m_2\frac{-C_1}{4B_1})}\\
\times\int_{N/(\Gamma_{\infty}\cap w\Gamma_{\infty}^{op}w^{-1})}
f_{\gamma^\prime}\left(n(0,y,z)w\right)\psi^{-1}\left(n(x,y,z)\right)dn.\label{SL3Rwhit12eqn6}
\end{multline}

If $m_1\neq0$, then this expression is 0, as can be seen by integrating over $x$. Thus suppose $m_1=0$. In this case,  
\[\int_{\Gamma_{\infty}\backslash N}\sum_{\gamma}\langle f,\pi((\gamma n)^{-1})\tau\rangle_{\lambda,\phi}\psi(n)dn\]
\[=\sum_{\gamma^\prime}|4B_1|^{-1-\lambda_1+\lambda_2}|4A_2|^{-1-\lambda_2+\lambda_3}\epsilon_{\gamma^\prime}s(\gamma^\prime)e^{2\pi i(m_2\frac{-C_1}{4B_1})}
\int_{\mathbb{R}^2}
 f_{\gamma^\prime}\left(n(0,y,z)
w\right)e^{2\pi i(-m_2y)}dydz.\]

Let $\varepsilon_{1},\varepsilon_{2}\in \{\pm1\}$ and define 
\begin{equation*}
S_{\varepsilon_{1},\varepsilon_{2}}=\sum_{\begin{smallmatrix}
\varepsilon_{1}B_1\in\mathbb{Z}_{>0}\\
\varepsilon_{2}A_2\in\mathbb{Z}_{>0}\\
4B_1|A_2
\end{smallmatrix}}
|4B_1|^{-1-\lambda_1+\lambda_2}|4A_2|^{-1-\lambda_2+\lambda_3}\sum_{\begin{smallmatrix}
C_1\,(\text{mod }4B_1)\\
C_2\,(\text{mod }4A_2)\\
C_j\equiv -1(\text{mod }4)
\end{smallmatrix}}\left(\frac{A_2/B_1}{-C_2}\right)\left(\frac{-B_1}{-C_1}\right)e^{2\pi im_2\frac{-C_1}{4B_1}}.
\end{equation*} We have just shown that 
\begin{align*}\langle f(\tilde{g}),\int_{\Gamma_{\infty}\backslash N}\tilde{E}_{w_{\alpha_1}w_{\alpha_2}}(n\tilde{g})\psi(n)&dn\rangle_{\lambda,\phi}\\
&=(S_{1,1}+(-i)S_{1,-1})\int_{\mathbb{R}^2}f_2\left(n(0,y,z)w\right)e^{-2\pi im_2y}dydz\\
&\phantom{=}+(S_{-1,-1}+(-i)S_{-1,1})\int_{\mathbb{R}^2}f_1\left(n(0,y,z)w\right)e^{-2\pi im_2y}dydz.
\end{align*}

We will simplify the Dirichlet series $S_{-1,-1}$; the other series are similar. We begin by considering the change of variable $B_1\mapsto -B_1$ to get 

\begin{equation}\label{alpha12sum2}S_{-1,-1} =\sum_{\begin{smallmatrix}
B_1\in\mathbb{Z}_{>0}\\
A_2\in\mathbb{Z}_{<0}\\
4B_1|A_2
\end{smallmatrix}}
|4B_1|^{-1-\lambda_1+\lambda_2}|4A_2|^{-1-\lambda_2+\lambda_3}\sum_{\begin{smallmatrix}
C_1\,(\text{mod }4B_1)\\
C_2\,(\text{mod }4A_2)\\
C_j\equiv -1(\text{mod }4)
\end{smallmatrix}}\left(\frac{-A_2/B_1}{-C_2}\right)\left(\frac{B_1}{C_1}\right)e^{2\pi im_2\frac{C_1}{4B_1}}.
\end{equation}
This is followed by the change of variables $k=\frac{-A_2}{4B_1}$.
\begin{equation}
(\ref{alpha12sum2}) = \sum_{\begin{smallmatrix}
B_1\in\mathbb{Z}_{>0}\\
k\in\mathbb{Z}_{>0}
\end{smallmatrix}}
|4B_1|^{-2-\lambda_1+\lambda_3}|4k|^{-1-\lambda_2+\lambda_3}\sum_{\begin{smallmatrix}
C_1\,(\text{mod }4B_1)\\
C_2\,(\text{mod }16kB_1)\\
C_j\equiv -1(\text{mod }4)
\end{smallmatrix}}\left(\frac{4k}{-C_2}\right)\left(\frac{B_1}{C_1}\right)e^{2\pi im_2\frac{C_1}{4B_1}}.
\end{equation}
Since the character $\left(\frac{4k}{-C_2}\right)$ only depends on $C_2$  modulo $4k$ we have 
\begin{equation}\label{alpha12sum3}
(\ref{alpha12sum2})=\sum_{\begin{smallmatrix}
B_1\in\mathbb{Z}_{>0}\\
k\in\mathbb{Z}_{>0}
\end{smallmatrix}}
|4B_1|^{-1-\lambda_1+\lambda_3}|4k|^{-1-\lambda_2+\lambda_3}\sum_{\begin{smallmatrix}
C_1\,(\text{mod }4B_1)\\
C_2\,(\text{mod }4k)\\
C_j\equiv -1(\text{mod }4)
\end{smallmatrix}}\left(\frac{4k}{-C_2}\right)\left(\frac{B_1}{C_1}\right)e^{2\pi im_2\frac{C_1}{4B_1}}.
\end{equation}
Now the sum can be factored.

\begin{multline}\label{alpha12sum4}
(\ref{alpha12sum3})=\left(
\sum_{B_1\in\mathbb{Z}_{>0}}
|4B_1|^{-1-\lambda_1+\lambda_3}\sum_{\begin{smallmatrix}
C_1\,(\text{mod }4B_1)\\
C_1\equiv -1(\text{mod }4)
\end{smallmatrix}}\left(\frac{B_1}{C_1}\right)e^{2\pi im_2\frac{C_1}{4B_1}}\right)\\
\times\left(\sum_{k\in\mathbb{Z}_{>0}}
|4k|^{-1-\lambda_2+\lambda_3}\sum_{\begin{smallmatrix}
C_2\,(\text{mod }4k)\\
C_2\equiv -1(\text{mod }4)
\end{smallmatrix}}\left(\frac{4k}{-C_2}\right)
\right).
\end{multline}
The first factor in line (\ref{alpha12sum4})  can be rewritten as  

\begin{multline}\label{alpha12sumfactor1}
\sum_{B_1\in\mathbb{Z}_{>0}}
|4B_1|^{-1-\lambda_1+\lambda_3}\sum_{\begin{smallmatrix}
C_1\,(\text{mod }4B_1)\\
C_1\equiv -1(\text{mod }4)
\end{smallmatrix}}\left(\frac{B_1}{C_1}\right)e^{2\pi im_2\frac{C_1}{4B_1}}\\
=i\sum_{B_1\in\mathbb{Z}_{>0}}
|4B_1|^{-1-\lambda_1+\lambda_3}\left(\frac{K_{1}(m_2;4B_1)-K_{-1}(m_2;4B_1)}{2}\right).
\end{multline}
As for the second factor in line (\ref{alpha12sum4}), the sum of quadratic characters will be nonzero precisely when $k$ is a square. Thus by the identity in line (\ref{squarephi})

\begin{equation}
\label{alpha12sumfactor2}
\sum_{k\in\mathbb{Z}_{>0}}
|4k|^{-1-\lambda_2+\lambda_3}\sum_{\begin{smallmatrix}
C_2\,(\text{mod }4k)\\
C_2\equiv -1(\text{mod }4)
\end{smallmatrix}}\left(\frac{4k}{-C_2}\right)
=2^{-2(-1-\lambda_2+\lambda_3)}\frac{\zeta(2\lambda_2-2\lambda_3)}{\zeta_2(2\lambda_2-2\lambda_3+1)}.
\end{equation}
By combining lines (\ref{alpha12sumfactor1}) and (\ref{alpha12sumfactor2}),  we see that 

\begin{multline}
S_{-1,-1}=i\left(\sum_{B_1\in\mathbb{Z}_{>0}}
|4B_1|^{-1-\lambda_1+\lambda_3}\left(\frac{K_{1}(m_2;4B_1)-K_{-1}(m_2;4B_1)}{2}\right)\right)\\
\times2^{-2(-1-\lambda_2+\lambda_3)}\frac{\zeta(2\lambda_2-2\lambda_3)}{\zeta_2(2\lambda_2-2\lambda_3+1)}.
\end{multline}
\EndProof

\subsection{Bruhat Cell: $Nw_{\alpha_{2}}w_{\alpha_{1}}B$} \label{ssec:tildeNw2w1B}


Let $w=\left(\left( \begin{smallmatrix}
1 & 0 & 0\\
0 & 0 & 1\\
0 & -1 & 0 \end{smallmatrix} \right), 1\right)
=\left(\left( \begin{smallmatrix}
 & &1\\
1 &  & \\
 & 1 &  \end{smallmatrix} \right)
w_{\ell}, 1\right)$. We begin by computing the effect of each summand of $\widetilde{E}_{w_{\alpha_2}w_{\alpha_1}}$ on a test vector.

\begin{proposition}\label{distmetaalpha2alpha1} Let $\gamma\in Nw_{\alpha_2}w_{\alpha_1}B$ with Pl\"{u}cker coordinates $(4A_1,4B_1,C_1,0,4B_2,C_2)$ and $f=
\left[\begin{smallmatrix}
f_1\\
f_2
\end{smallmatrix}\right]\in \widetilde{V}_{-\lambda,^\intercal\phi^{-1}}^{\infty}$. Then:
\begin{multline*}
\langle f,\pi(\gamma^{-1})\tilde{\tau}\rangle_{\lambda,\phi}=|4A_1|^{-1-\lambda_1+\lambda_2}|4B_2|^{-1-\lambda_2+\lambda_3}\left(\frac{-A_1/B_2}{-C_1}\right)\left(\frac{B_2}{-C_2}\right)\\
\times\begin{cases}
\phantom{-i}f_2\left(wn(\frac{C_1}{4A_1},0, \frac{-B_1}{A_1})\right),&\text{ if } A_1>0,B_2>0;\\
\phantom{-}if_2\left(wn(\frac{C_1}{4A_1},0, \frac{-B_1}{A_1})\right),&\text{ if } A_1>0,B_2<0;\\
-\phantom{i}f_1\left(wn(\frac{C_1}{4A_1},0, \frac{-B_1}{A_1})\right),&\text{ if } A_1<0,B_2<0;\\
\phantom{-}if_1\left(wn(\frac{C_1}{4A_1},0, \frac{-B_1}{A_1})\right),&\text{ if } A_1<0,B_2>0.
\end{cases}
\end{multline*}
\end{proposition}

Now we compute the Fourier coefficients of $\widetilde{E}_{w_{\alpha_2}w_{\alpha_1}}$.


 \begin{proposition}\label{metawhitalpha21}
Let $f=
\left[\begin{smallmatrix}
f_1\\
f_2
\end{smallmatrix}\right]\in \widetilde{V}_{-\lambda,^\intercal\phi^{-1}}^{\infty}$.

If $m_2\neq0$, then $\langle f(\tilde{g}),\int_{\Gamma_{\infty}\backslash N}\tilde{E}_{w_{\alpha_2}w_{\alpha_1}}(n\tilde{g})\psi(n)dn\rangle_{\lambda,\phi}=0.$

If $m_2=0$, then
\begin{multline*} 
\langle f(\tilde{g}),\int_{\Gamma_{\infty}\backslash N}\tilde{E}_{w_{\alpha_2}w_{\alpha_1}}(n\tilde{g})\psi(n)dn\rangle_{\lambda,\phi}
=\\
\left(2^{-2(-1-\lambda_1+\lambda_2)}\frac{\zeta(2\lambda_1-2\lambda_2)}{\zeta_2(2\lambda_1-2\lambda_2+1)}\right)
\Bigg(\sum_{
B_2\in\mathbb{Z}_{>0}}|4B_2|^{-1-\lambda_1+\lambda_3}
K_{1}(m_1;4B_2)\Bigg)\\
\times\int_{\mathbb{R}^2}(-f_1+if_2)\left(n(x,0,z)
\left(\left(\begin{smallmatrix}
1 & 0 & 0\\
0 & 0 & 1\\
0 & -1 & 0 \end{smallmatrix}\right),1\right)\right)\psi^{-1}(n(x,0,z))dxdz.
\end{multline*}
\end{proposition}


\subsection{Bruhat Cell: $Nw_{\ell}B$} \label{ssec:tildeNwlB}


Let $w=\left(\left(\begin{smallmatrix}
-1 & & \\
 & 1 & \\
 & & -1 \end{smallmatrix}\right),1\right)=\left(\left( \begin{smallmatrix}
  &  & -1\\
 & -1 & \\
-1 &  &  \end{smallmatrix} \right)w_{\ell},1\right)$. We begin by computing the effect of each summand of $\widetilde{E}_{w_{\ell}}$ on a test vector.

\begin{proposition}\label{distmetabigcell} If $\gamma\in Nw_{\ell}B$ with Pl\"{u}cker coordinates $(4A_1,4B_1,C_1,4A_2,4B_2,C_2)$ and $f=
\left[\begin{smallmatrix}
f_1\\
f_2
\end{smallmatrix}\right]\in \widetilde{V}_{-\lambda,^\intercal\phi^{-1}}^{\infty}$, then:
\begin{multline*}
\langle f,\pi(\gamma^{-1})\tilde{\tau}\rangle_{\lambda,\phi}=|4A_1|^{-1-\lambda_1+\lambda_2}|4A_2|^{-1-\lambda_2+\lambda_3}s(\gamma)\\
\times\begin{cases}
\phantom{-i}f_1\left(wn(\frac{B_1}{A_1},\frac{-B_2}{A_2},\frac{C_2}{4A_2})\right),&\text{ if } A_1>0,A_2>0;\\
\phantom{-}if_1\left(wn(\frac{B_1}{A_1},\frac{-B_2}{A_2},\frac{C_2}{4A_2})\right),&\text{ if } A_1>0,A_2<0;\\
-\phantom{i}f_2\left(wn(\frac{B_1}{A_1},\frac{-B_2}{A_2},\frac{C_2}{4A_2})\right),&\text{ if } A_1<0,A_2<0;\\
-if_2\left(wn(\frac{B_1}{A_1},\frac{-B_2}{A_2},\frac{C_2}{4A_2})\right),&\text{ if } A_1<0,A_2>0.
\end{cases}
\end{multline*}

Recall that the formula for $s(\gamma)$ in this case is contained in Theorem \ref{theorem:PluckSplit}.

\end{proposition}

Now we compute the Fourier coefficients of $\widetilde{E}_{w_{\ell}}$.

\begin{proposition}\label{metawhitwell}
Let $f=
\left[\begin{smallmatrix}
f_1\\
f_2
\end{smallmatrix}\right]\in \widetilde{V}_{-\lambda,^\intercal\phi^{-1}}^{\infty}$.\newline
\begin{align*} 
&\hspace{5cm}\langle f(\tilde{g}),\int_{\Gamma_{\infty}\backslash N}\tilde{E}_{w_{\ell}}(n\tilde{g})\psi(n)dn\rangle_{\lambda,\phi}\\
=&\sum_{\begin{smallmatrix}
A_1>0\\
A_2>0
\end{smallmatrix}}\
|4A_1|^{-1-\lambda_1+\lambda_2}|4A_2|^{-1-\lambda_2+\lambda_3}(\Sigma(A_1,A_2;-m_1,m_2)+i\Sigma(A_1,-A_2;-m_1,m_2))\\
&\hspace{5cm}\times
\int_{\mathbb{R}^3}f_1\left(n(x,y,z)
w\right)\psi^{-1}(n(x,y,z))dxdydz\\
&+\sum_{\begin{smallmatrix}
A_1>0\\
A_2>0
\end{smallmatrix}}|4A_1|^{-1-\lambda_1+\lambda_2}|4A_2|^{-1-\lambda_2+\lambda_3}(-\Sigma(A_1,A_2;-m_1,-m_2)-i\Sigma(A_1,-A_2;-m_1,-m_2))\\
&\hspace{5cm}\times
\int_{\mathbb{R}^3}f_2\left(n(x,y,z)
w\right)\psi^{-1}(n(x,y,z))dxdydz.
\end{align*}
\end{proposition}

The terms $\Sigma(A_1,A_2;m_1,m_2)$ satisfy a twisted multiplicativity in $A_1$ and $A_2$, stated in Proposition \ref{exptwistmult}; a form of twisted multiplicativity in $m_1$ and $m_2$, stated in Proposition \ref{twistmultindex}; and the symmetries of Proposition \ref{ExpSym}. Thus the computation of $\Sigma(A_1,A_2;m_1,m_2)$ may be reduced to that of $\Sigma(p^k,\mu p^l;p^{r_1},p^{r_2})$, where $\mu=\pm1$. Formulas for these expressions may be found in Section \ref{sec:BBFHform}.

\noindent\textbf{Proof:} In what follows, the summation over $\gamma$ will be described by P\"{u}cker coordinates for the set $\Gamma_{\infty}\backslash \Gamma\cap(Nw_{\ell}B)$. The summation over $\gamma^\prime$ will consist of elements of the double coset space $\Gamma_{\infty}\backslash \Gamma\cap(Nw_{\ell}B)/\Gamma_{\infty}$. 

Let $n$ denote $n(x,y,z)$. Define $\epsilon_{\gamma}\in\{\pm1,\pm i\}$ and $f_{\gamma}=f_{j}$, where $j\in\{1,2\}$ be defined so that 
\[\langle \pi(n)f,\pi(S(\gamma)^{-1})\tilde{\tau}\rangle_{\lambda,\phi}=|4A_1|^{-1-\lambda_1+\lambda_2}|4A_2|^{-1-\lambda_2+\lambda_3}\epsilon_\gamma s(\gamma)
f_{\gamma}\left(wn(\frac{B_1}{A_1},\frac{-B_2}{A_2},\frac{C_2}{4A_2})\right),\]
in accordance with Proposition \ref{distmetabigcell}. Observe that $\epsilon_{\gamma}$, $s(\gamma)$, and $j$ only depend on the double coset $\Gamma_{\infty}\backslash \Gamma /\Gamma_{\infty}$. In fact $\epsilon_{\gamma}$ and $f_{\gamma}$ only depend on the signs of $A_1$ and $A_2$, so we will write $\epsilon_{A_1,A_2}=\epsilon_{\gamma}$ and $f_{A_1,A_2}=f_{\gamma}$.

Begin with the change of variables $n\mapsto n^{-1}$ and then apply Proposition \ref{distmetabigcell} to see that
\begin{multline}
\langle f,\int_{\Gamma_{\infty}\backslash N}\tilde{E}_{w_{\ell}}(ng)\,\psi(n)dn\rangle_{\lambda,\phi}\\
=\int_{N/\Gamma_{\infty}}\sum_{\gamma}|4A_1|^{-1-\lambda_1+\lambda_2}|4A_2 |^{-1-\lambda_2+\lambda_3}\epsilon_{\gamma}s(\gamma)f_{\gamma}\left(nwn(\frac{B_1}{A_1}, \frac{-B_2}{A_2},\frac{C_2}{4A_2})\right)\psi^{-1}(n)dn.\label{SL3Rwhitelleqn1}
\end{multline}
Since the 2-cocycle $\sigma$ is trivial on $N$ we have $wn(\frac{B_1}{A_1}, \frac{-B_2}{A_2},\frac{C_2}{4A_2})=n(\frac{-B_1}{A_1}, \frac{B_2}{A_2},\frac{C_2}{4A_2})w$. Thus 
\begin{equation}(\ref{SL3Rwhitelleqn1})=\int_{N/\Gamma_{\infty}}\sum_{\gamma}|4A_1|^{-1-\lambda_1+\lambda_2}|4A_2|^{-1-\lambda_2+\lambda_3}\epsilon_{\gamma}s(\gamma)f_{\gamma}\left(nn(\frac{-B_1}{A_1}, \frac{B_2}{A_2},\frac{C_2}{4A_2})w\right)\psi^{-1}(n)dn.\label{SL3Rwhitelleqn2}\end{equation}
The next step is to unfold the integral. For details recall the analogous step in Proposition \ref{metawhitalpha12}. In this case,
\begin{equation}(\ref{SL3Rwhitelleqn2})=\sum_{\gamma^\prime}|4A_1|^{-1-\lambda_1+\lambda_2}|4A_2 |^{-1-\lambda_2+\lambda_3}\epsilon_{\gamma^\prime}s(\gamma^\prime)\int_{N}f_{\gamma}\left(nn(\frac{-B_1}{A_1}, \frac{B_2}{A_2},\frac{C_2}{4A_2})w\right)\psi^{-1}(n)dn.\label{SL3Rwhitelleqn3}\end{equation}
Finally, apply the change of variables $n\mapsto nn(\frac{-B_1}{A_1}, \frac{B_2}{A_2},\frac{C_2}{4A_2})^{-1}$ to get
\begin{equation}(\ref{SL3Rwhitelleqn3})=\sum_{\gamma^\prime}|4A_1|^{-1-\lambda_1+\lambda_2}|4A_2|^{-1-\lambda_2+\lambda_3}\epsilon_{\gamma^\prime}s(\gamma^\prime)e^{2\pi i(m_1\frac{-B_1}{A_1}+m_2\frac{B_2}{A_2})}\int_{N}f_{\gamma^\prime}\left(nw\right)\psi^{-1}(n)dn.\label{SL3Rwhitelleqn4}\end{equation}

Now line $(\ref{SL3Rwhitelleqn4})$ is equal to
\begin{equation}
\sum_{A_1,A_2\in\mathbb{Z}_{\neq 0}}|4A_1|^{-1-\lambda_1+\lambda_2}|4A_2|^{-1-\lambda_2+\lambda_3}\epsilon_{A_1,A_2}\Sigma(A_1,A_2;-m_1,m_2)\int_{N}f_{A_{1},A_2}\left(nw\right)\psi^{-1}(n)dn.
\end{equation}
We can complete the proof by applying Proposition \ref{ExpSym} and by using Proposition \ref{distmetabigcell} to evaluate $\epsilon_\gamma$, and $f_{\gamma}$. \EndProof

\subsection{Constant Term} \label{ssec:Constant}

The computations of the previous section can be specialized ($m_1=m_2=0$) to produce the constant term. 

\begin{theorem}\label{theorem:Constant}
Let $f=
\left[\begin{smallmatrix}
f_1\\
f_2
\end{smallmatrix}\right]\in \widetilde{V}_{-\lambda,^\intercal\phi^{-1}}^{\infty}$.

\begin{align*} 
&\hspace{5cm}\langle f(\tilde{g}),\int_{\Gamma_{\infty}\backslash N}\tilde{E}(n\tilde{g})dn\rangle_{\lambda,\phi}=f_1((w_{\ell},1))\\
&+(1-i)
2^{-2(1+\lambda_2-\lambda_3)}\frac{\zeta(2(\lambda_2-\lambda_3))}{\zeta_{2}(2(\lambda_2-\lambda_3)+1)}
\int_{\mathbb{R}}f_2(n(x,0,0)\left(
\left(\begin{smallmatrix}
  & -1 & \\
-1 &  & \\
&  & -1 \end{smallmatrix} \right)
w_{\ell},1\right))dx\\
&+2^{-2(1+\lambda_1-\lambda_2)}\frac{\zeta(2(\lambda_1-\lambda_2))}{\zeta_{2}(2(\lambda_1-\lambda_2)+1)}
\int_{\mathbb{R}}(if_1-f_2)\left(n(0,y,0)
\left(
\left( \begin{smallmatrix}
-1 &  & \\
&  & -1\\
& -1 &  \end{smallmatrix} \right)
w_{\ell},1\right)
\right)dy\\
&+(1-i)2^{-2(1+\lambda_1-\lambda_3)}2^{-2(1+\lambda_2-\lambda_3)}\frac{\zeta(2(\lambda_1-\lambda_3))\zeta(2\lambda_2-2\lambda_3)}{\zeta_{2}(2(\lambda_1-\lambda_3)+1)\zeta_2(2(\lambda_2-\lambda_3)+1)}\\
&\hspace{7cm}\times\int_{\mathbb{R}^2}(f_1+if_2)\left(n(0,y,z)\left(
\left( \begin{smallmatrix}
 & 1 & \\
&  & 1\\
1&  &  \end{smallmatrix} \right)
w_{\ell},1\right)\right)dydz\\
&+(1-i)2^{-2(1+\lambda_1-\lambda_3)}2^{-2(1+\lambda_1-\lambda_2)}\frac{\zeta(2(\lambda_1-\lambda_3))\zeta(2(\lambda_1-\lambda_2))}{\zeta_2(2(\lambda_1-\lambda_3)+1)\zeta_2(2(\lambda_1-\lambda_2)+1)}\\
&\hspace{7cm}\times\int_{\mathbb{R}^2}(-f_1+if_2)\left(n(x,0,z)
\left(\left( \begin{smallmatrix}
 & & 1\\
1 &  & \\
 & 1 &  \end{smallmatrix} \right)
w_{\ell},1\right)\right)dxdz\\
&+i2^{-2(1+\lambda_1-\lambda_2)}2^{-2(1+\lambda_2-\lambda_3)}(1-2^{-2(\lambda_1-\lambda_2)}-2^{-2(\lambda_2-\lambda_3)}+6(2^{-2(\lambda_1-\lambda_3+1)}))\\
&\hspace{4cm}\times\frac{\zeta(2(\lambda_1-\lambda_2))\zeta(2(\lambda_2-\lambda_3))\zeta(2(\lambda_1-\lambda_3))}{\zeta_2(2(\lambda_1-\lambda_2)+1)\zeta_2(2(\lambda_2-\lambda_3)+1)\zeta_2(2(\lambda_1-\lambda_3)+1)}\\
&\hspace{6cm}\times\int_{\mathbb{R}^3}(f_1-f_2)\left(n(x,y,z)
\left(\left(\begin{smallmatrix}
 &  & -1 \\
 & -1 & \\
-1 &  & \end{smallmatrix}\right)w_{\ell},1\right)\right)dxdydz.
\end{align*}

This formula can also be written more succinctly as follows: Let
\begin{equation*}
F_{\ell,2}(\lambda_1,\lambda_2,\lambda_3) = 2^{2+2(\lambda_1-\lambda_3)}(1-2^{-2(\lambda_1-\lambda_2)}-2^{-2(\lambda_2-\lambda_3)}+6(2^{-2(\lambda_1-\lambda_3+1)})),
\end{equation*} 
\[\begin{array}{rclcrcll}
v_{id}&=&
^\intercal\left[\begin{smallmatrix}
1\\
0
\end{smallmatrix}\right]&,
&v_{w_{\ell}}&=&
F_{\ell,2}(\lambda)^\intercal\left[\begin{smallmatrix}
i\\
-i
\end{smallmatrix}\right]&,\\
v_{w_{\alpha_1}}&=&
^\intercal\left[\begin{smallmatrix}
0\\
1-i
\end{smallmatrix}\right]&,
&v_{w_{\alpha_2}}&=&
^\intercal\left[\begin{smallmatrix}
i\\
-1
\end{smallmatrix}\right]&,\\
v_{w_{\alpha_1}w_{\alpha_2}}&=&
^\intercal\left[\begin{smallmatrix}
1-i\\
1+i
\end{smallmatrix}\right]&,
&v_{w_{\alpha_2}w_{\alpha_1}}&=&
^\intercal\left[\begin{smallmatrix}
-1+i\\
1+i
\end{smallmatrix}\right]&,
\end{array}\] and $W^\prime=\left\{
\left(\begin{smallmatrix}
1& &\\
& 1&\\
& &1
\end{smallmatrix}\right),
\left(\begin{smallmatrix}
& -1&\\
-1& &\\
& &1
\end{smallmatrix}\right),
\left(\begin{smallmatrix}
1& &\\
& &-1\\
& -1&
\end{smallmatrix}\right),
\left(\begin{smallmatrix}
& & 1\\
1& &\\
& 1&
\end{smallmatrix}\right),
\left(\begin{smallmatrix}
& 1&\\
& &1\\
1& &
\end{smallmatrix}\right),
\left(\begin{smallmatrix}
& &-1\\
& -1&\\
-1& &
\end{smallmatrix}\right)\right\}$. Then:

\begin{align*}
&\hspace{5cm}\langle f(\tilde{g}),\int_{\Gamma_{\infty}\backslash N}\tilde{E}(n\tilde{g})dn\rangle_{\lambda,\phi}\\
&=\sum_{w\in W^\prime}\left(\prod_{\alpha\in \Phi^{+}\cap ww_{\ell}\Phi^{+}}4^{-1-(\lambda,h_\alpha)}\frac{\zeta(2(\lambda,h_{\alpha}))}{\zeta_{2}(2(\lambda,h_{\alpha})+1)}\right)\\
&\hspace{5cm}\times\int_{w^{-1}w_{\ell}N(w^{-1}w_{\ell})^{-1}\cap N}\left(
v_{w}\cdot
\left[\begin{smallmatrix}
f_1\\
f_2
\end{smallmatrix}\right]\right)((n,1)(w^{-1}w_{\ell},1))dn.
\end{align*}

\end{theorem}

\textbf{Proof:} To compute the portion of the constant term that does not arise from the big cell $Nw_{\ell}B$ we may set $m_1=m_2=0$ in propositions \ref{metawhitid}, \ref{metawhitalpha1}, \ref{metawhitalpha2}, \ref{metawhitalpha12}, and \ref{metawhitalpha21} and apply the identities of lines (\ref{squarephi}) and (\ref{KloostermanConstant}). 

Finally, we consider the big cell $Nw_{\ell}B$. We begin by setting $m_1=m_2=0$ in Proposition \ref{metawhitwell}. From Proposition \ref{ExpSym} we see that $\Sigma(A_1,A_2;m_1,m_2)=0$ if $A_1A_2>0$. Thus we may focus on computing $\Sigma(A_1,-A_2;0,0)$ where $A_1,A_2>0$. Note that Proposition \ref{ExpSym} also implies that $\Sigma(A_1,-A_2;0,0)=\Sigma(A_2,-A_1;0,0)$. By twisted multiplicativity, Proposition \ref{twistmult}, we see that up to sign $\Sigma(A_1,-A_2;0,0)$ is a product of $\Sigma(p^k,\pm p^{\ell};0,0)$ for $p$ prime. The computations of Subsection \ref{ssec:bigcellConstant} show us that if $\Sigma(p^k,\epsilon p^{\ell};0,0)\neq0$, then $k$ and $\ell$ are even and $\epsilon=-1$. Thus the twisted multiplicativity of Proposition \ref{twistmult} is a true mutliplicativity and 

\begin{equation}\sum_{\begin{smallmatrix}
A_1>0\\
A_2>0
\end{smallmatrix}}\
|4A_1|^{-1-\lambda_1+\lambda_2}|4A_2|^{-1-\lambda_2+\lambda_3}\Sigma(A_1,-A_2;0,0)
\end{equation}
has an Euler product. Thus we may focus on computing 

\begin{equation}\label{constppart}\sum_{\begin{smallmatrix}
k>0\\
\ell>0
\end{smallmatrix}}\
|p^{2k}|^{-1-\lambda_1+\lambda_2}|p^{2\ell}|^{-1-\lambda_2+\lambda_3}\Sigma(p^{2k},-p^{2\ell};0,0),
\end{equation}
where $p$ is a prime. We can write line \ref{constppart} as a rational function using the computations of Subsection \ref{ssec:bigcellConstant}. First we will consider the case where $p$ is an odd prime.

In this case we can evaluate $\Sigma(p^{2k},-p^{2\ell};0,0)$ using propositions \ref{sigma+-p k,0,0,0}, \ref{sigma+-p k>l>0,0,0}, \ref{sigma+-p,k>0,0,0}, and \ref{sigma++,p,k>0,0,0}. After applying the formula for geometric series and some algebraic simplifications we find that when $p$ is an odd prime

\begin{equation}
(\ref{constppart})=\frac{(1-p^{-2(\lambda_1-\lambda_2)-1})(1-p^{-2(\lambda_2-\lambda_3)-1})(1-p^{-2(\lambda_1-\lambda_3)-1})}{(1-p^{-2(\lambda_1-\lambda_2)})(1-p^{-2(\lambda_2-\lambda_3)})(1-p^{-2(\lambda_1-\lambda_3)})}.
\end{equation}

We proceed in the same manner for the prime $p=2$. In this case we can evaluate $\Sigma(2^{2k},-2^{2\ell};0,0)$ using propositions \ref{sigma 2,k>l>0,0,0}, and \ref{sigma 2,k>0,0,0}. When $p=2$ we find that 

\begin{equation}
(\ref{constppart})=\frac{1-2^{-2(\lambda_1-\lambda_2)}-2^{-2(\lambda_2-\lambda_3)}+6(2^{-2(\lambda_1-\lambda_3+1)})}{(1-2^{-2(\lambda_1-\lambda_2)})(1-2^{-2(\lambda_2-\lambda_3)})(1-2^{-2(\lambda_1-\lambda_3)})}.
\end{equation}

The result follows once we identify the Euler product as a product of zeta functions.\EndProof

\section{Acknowledgements}
I would like to thank Ben Brubaker, Daniel Bump, Gautam Chinta, Solomon Friedberg, Paul Gunnells, Jeff Hoffstein, Henryk Iwaniec, and Martin Weissman for many valuable conversations during the preparation of this paper. I thank the referee for many helpful comments and suggestions. I would also like to thank Stephen D. Miller for his guidance and for access to his unpublished notes, which provided the starting point for this project.

Declarations of interest: none

\end{document}